%% file: L2HodgeDeRham.tex
\documentclass[11pt]{article}
\usepackage{amsmath,amsthm,amsfonts,latexsym,amscd}
\renewcommand{\qed}{q.e.d.}

\usepackage{amsthm}

\swapnumbers
\theoremstyle{plain}
\newtheorem{theorem}{Theorem}[section]
\newtheorem{lemma}[theorem]{Lemma}
\newtheorem{corollary}[theorem]{Corollary}
\newtheorem{proposition}[theorem]{Proposition}

\theoremstyle{definition}
\newtheorem{definition}[theorem]{Definition}
\newtheorem{example}[theorem]{Example}
\newtheorem{notation}[theorem]{Notation}
\newtheorem{convention}[theorem]{Convention}

\theoremstyle{remark}
\newtheorem{remark}[theorem]{Remark}

\usepackage{amsmath,amsfonts,latexsym}

\newcommand{\R}{\mathbb{R}}
\newcommand{\reals}{\mathbb{R}}
\newcommand{\C}{\mathbb{C}}

\newcommand{\N}{\mathbb{N}}
\newcommand{\naturals}{\mathbb{N}}
\newcommand{\Z}{\mathbb{Z}}
\newcommand{\integers}{\mathbb{Z}}

\newcommand{\K}{\mathbb{K}}

\newcommand{\1}{\mathbf{1}}  

\newcommand{\boundary}[1]{\partial#1}
\newcommand{\boundedops}{\mathcal{B}}
\newcommand{\bundleover}{\!\!\downarrow\!\!}
\newcommand{\abs}[1]{\left\lvert#1\right\rvert} 
\newcommand{\norm}[1]{\left\lVert#1\right\rVert}

\newcommand{\tensor}{\otimes}

\DeclareMathOperator{\supp}{supp}   
\DeclareMathOperator{\im}{im}      
\DeclareMathOperator{\vol}{vol}    
\DeclareMathOperator{\dist}{dist}  
\DeclareMathOperator{\End}{End}    
\DeclareMathOperator{\Hom}{Hom}    

\DeclareMathOperator{\trace}{trace}
\DeclareMathOperator{\tr}{tr}

\DeclareMathOperator{\Ric}{Ric}  

\newcommand{\forget}[1]{}

{\catcode`@=11\global\let\c@equation=\c@theorem}

\renewcommand{\labelenumi}{(\arabic{enumi})}

\allowdisplaybreaks[2]

\newcommand{\bwp}[1]{\mathcal{#1}}
\newcommand{\bvp}[1]{\mathcal{#1}}
\newcommand{\vect}[1]{\vec{#1}}
\newcommand{\domain}{\mathcal}

\newcommand{\bundlesum}[1]{\mathcal{#1}}

\begin{document}
\input{L2_intro}
\input{bounded_geo}
\input{elliptic}
\input{laplace}
\input{L2+curvature}
\input{L2_DodziukSchick}

\bibliographystyle{lueck}
\bibliography{literatur}
\end{document}

%% file: L2_intro.tex
\title{Geometry and Analysis of Boundary-Manifolds of Bounded Geometry}
\author{Thomas Schick\thanks{e-mail:
    thomas.schick@math.psu.edu}\\Dept. of Mathematics --- Penn State
University\\
218 McAllister Building ---
         University Park, PA 16802} 
         
\maketitle

\begin{abstract} In this paper, we investigate analytical and
geometric properties of certain non-compact boundary-manifolds, namely
manifolds of bounded geometry.  One result are strong Bochner type
vanishing results for the $L^2$-coho\-mo\-lo\-gy of these manifolds:
if e.g.~a
manifold admits a metric of bounded geometry which outside a compact
set has nonnegative Ricci curvature and nonnegative mean curvature (of
the boundary) then its first relative $L^2$-coho\-mo\-logy vanishes (this
in particular answers a question of Roe).

We prove the Hodge-de Rham-theorem for $L^2$-cohomology of oriented
$\boundary$-manifolds of bounded geometry.

The technical basis is the study of (uniformly elliptic) boundary
value problems on these manifolds, applied to the Laplacian.

\end{abstract}


\section{Introduction}
The broad theme of this paper is geometry, topology and analysis on
non-compact manifolds with boundary, which might be empty. We focus on
manifolds of bounded geometry (Definition \ref{def_bounded}). These
are (non-compact) Riemannian manifolds with bounds on curvature and
second fundamental form and their derivatives. Also, they have
positive injectivity radius.

Note that, on a non-compact manifold, very different metrics may
exist. Therefore, we will concentrate on obstructions for Riemannian
metrics in a given bilipschitz class. This is proposed by Roe in
\cite{Roe(1988c)}.

We use the Weizenb\"ock formula (for boundary-manifolds) to derive the
following Bochner vanishing Theorem \ref{Boch}: if we have a metric
with nonnegative Ricci curvature and convex boundary, then the first
$L^2$-cohomology is trivial. Then we refine this method so that we can
weaken the positivity conditions: they are not required to hold on all
of the manifold but only on an eventually large set, defined as
follows:

\begin{definition}\label{eventlarge}
Let $M$ be any metric space. A subset $X\subset M$ is called {\em
eventually large} if  $\forall R>0$ we find $x_R\in X$ so
that $B(x_R,R)\subset X$; the $R$-ball around $x_R$ is entirely
contained in $X$.
\end{definition}
 
A corollary is 
the following generalization of results of Roe \cite{Roe(1988c)}:
\begin{theorem} (compare Theorem \ref{str_Boch})\\
Let $M$ be a compact manifold, possibly with boundary. Let $\bar M$ be a normal
covering of $M$. Suppose we find a metric $g$ on $\bar M$
which lies in the bilipschitz class of the lifted metric and has
bounded geometry. Suppose there exists an eventually large subset
$X\subset \bar M$ so that on $X$ $g$ has nonnegative Ricci curvature and
nonnegative mean curvature (of the boundary). Then $b^1_{(2)}(\bar
M,\boundary \bar M)=0$.
\end{theorem}

\begin{corollary}
If $\dim M=2$, or $\dim M=4$ and $M$ is closed, and if
$\chi(M)<0$ then  on the universal
covering of $M$ no metric $g$ as described in the theorem exists.
\end{corollary}

These results answer affirmative  \cite[Question 6.3]{Roe(1988c)} of Roe.

The first (relative) $L^2$-Betti number $b_{(2)}^1$ of the theorem above  is zero exactly if the space of
$L^2$-harmonic $1$-forms (with relative boundary conditions) is trivial.
 It is the von Neumann dimension of this space  and a homotopy invariant of $M$ (compare \cite{Dodziuk(1977)}).
We prove for an arbitrary
orientable manifold $M$ of bounded geometry the $L^2$-Hodge de Rham
theorem:  the square integrable harmonic forms are isomorphic to
simplicially defined $L^2$-cohomology. 
For this end let $\boundary M=M_1\amalg M_2$ ($M_1$ or $M_2$ or both may be empty). 

Let $K$ be a smooth triangulation of $M$. It has to be sufficiently regular, adapted to the bounded geometry condition (compare \cite[2.3]{Dodziuk(1981)}):
\begin{definition}\label{triang}
The triangulation $K$ is called {\em $g$-bounded} if there are $\nu,K,c >0$ so that
\begin{enumerate}
\item the volume of the top-simplexes is bounded from below $v$,
\item the diameter of the top-simplexes is bounded from above by $K$,
\item $\abs{d\varphi_p}\le c$ $\forall p$ where $\varphi_p$ is the barycentric coordinate function for the vertex $p$ of $K$. 
\end{enumerate}
\end{definition} 
This  mimics the injectivity radius and curvature conditions
in the definition of bounded geometry.

\begin{definition}
The $L^2$-cochain complex $C^*_{(2)}(K,K_1)$ is defined as the subcomplex of the cochain complex $C^*(K,K_1)$ consisting of those chains whose coefficients are square summable (obviously a  Hilbert space) ($K_1$ is the  restriction of $K$ to $M_1$). Since $K$ is $g$-bounded this is  a sub-complex (compare \cite{Dodziuk(1981)}). Its $L^2$-cohomology is
\[ H^p_{(2)}(K,K_1) := \ker d_p/\overline{\im d_{p-1}} . \]
\end{definition}

The following is our  $L^2$-Hodge-de Rham theorem:

\begin{theorem}\label{rep_Rham_theorem}\label{Hodge_de_Rham_theorem}
Let ${\mathcal H^*}(M,  M_1)$ be the space of smooth harmonic square integrable forms on $M$ which fulfill absolute boundary conditions on $ M_2$ and relative boundary conditions on $ M_1$ (defined in \ref{def_of_harmonic}). Then the de Rham map (integration of forms over simplexes) 
\[A:\mathcal{H}^*( M, M_1)\to H^*_{(2)}(K, K_1)\]
 yields an isomorphism between this space and the relative $L^2$-cohomology $H^*_{(2)}(K,  K_1)$ of the triangulation $K$.\end{theorem}
If $\boundary
M=\emptyset$ this was done by Dodziuk \cite{Dodziuk(1981)}. We use
his result and a doubling trick to prove the result if the metric is
a product near the boundary. The last step consists in the comparison of an arbitrary metric with a suitable product metric near the boundary.

Basic technical tools to derive these results are ``uniform analysis''
like the introduction of uniform Sobolev spaces on manifolds of
bounded geometry  which have properties
similar to the Euclidian case. We introduce uniformly elliptic
boundary value problems for which we derive e.g.~regularity results as
usual. We prove that the Laplacian with absolute/relative boundary
conditions belongs to this class. This in turn implies that it is essentially self adjoint. 
The general theory also gives the following Hodge
decomposition (of Sobolev spaces):


\begin{theorem} (compare \ref{HsHodge} for the definitions)\\
Let $M$ be a manifold of bounded geometry, $\boundary M=M_1\amalg M_2$.
If $\mathcal H^p(M,M_1)$ is the space of square integrable harmonic forms which fulfill the Dirichlet/Neumann boundary conditions $\vect p$ then 
\[
H^{2k}(\Lambda^p(T^*M), bd(2k))={\mathcal H^p}(M, M_1) \oplus \overline{d\Omega^{p-1}_\infty} \oplus \overline{\delta\Omega^{p+1}_\infty} .
\]
The space which  we decompose consists of forms which fulfill Dirichlet/Neumann boundary conditions up to order $2k$. 
\end{theorem}

This paper is based on parts of the author's dissertation \cite{Schick(1996)}. I thank my advisor, Prof.~W.\ L\"uck, for his constant support and encouragement.

\subsection{Organization of the paper}
We start in Section \ref{sec_not} with some notational conventions which will be used throughout the text.

In Section \ref{sec_bound} we introduce the concept of a $\boundary$-manifold of bounded geometry and derive their main (analytical) properties.

Section \ref{sec_ell} deals with elliptic boundary value problems. We define uniform ellipticity as the right condition for manifolds of bounded geometry and derive basic properties.

In Section \ref{sec_lap} we proof the Hodge decomposition theorem.

After the technical tools are established,
Section \ref{sec_curv} examines the relations between $L^2$-cohomology
and curvature and proves the Bochner vanishing theorems.

In Section \ref{sec_Rham} we prove the $L^2$-de Rham theorem. 

\section{Notation}\label{sec_not}
\begin{definition}
We will often use inequalities to show equivalences of norms. If $a,b,c\in\R$ with $c>0$, set
$a\stackrel{c}{\le} b\quad:\iff a\le c\cdot b$.

If $c$ is not explicitly specified, a suitable positive constant has to be chosen. The same symbol (e.g.\ $c$) may be used for different constants.
\end{definition}

\begin{definition}
We will come up with equations where the sign of some of the terms does not matter. Then we simply use ``$\pm$''.
\end{definition}

\begin{convention}
A manifold $M$ shall always have dimension $m$, a bundle $E\bundleover M$ dimension $n$, if not stated otherwise. $\K$ stands for $\R$ or $\C$.
\end{convention}

\begin{definition}
We use multiindices $\alpha=(\alpha_1,\dots,\alpha_m)$: tupels of nonnegative integers. Set $\abs{\alpha}:=\max\{\alpha_1,\dots,\alpha_m\}$ and $D^\alpha:=\frac{\partial^{\alpha_1}}{\partial x_1^{\alpha_1}}\dots\frac{\partial^{\alpha_m}}{\partial x_m^{\alpha_m}}$.
\end{definition}

\begin{definition}
For a metric space $M$ set $B(p,r)=\{x\in M;\; d(x,p)\le r\}$. If we write $B(x',r)\subset \boundary M$ we consider $\boundary M$ as a metric space defined by the induced Riemannian metric and construct the ball in this metric.
\end{definition}

\begin{notation}
Let $M$ be a manifold of bounded geometry. A constant is called
$M$-universal if it depends only on $M$, not on local data.
\end{notation}

%% file: bounded_geo.tex
\section{Bounded geometry}\label{sec_bound}

In this section we introduce the concept of a boundary-manifold of
bounded geometry. These manifolds have uniformity properties which
e.g.~allow the definition of uniform Sobolev spaces. For them, the
Sobolev embedding theorem holds; differential operators induce bounded
operators $H^s\to H^{s-k}$; for $s>1/2$ a bounded restriction map
$H^s(M)\to H^{s-1/2}(\boundary M)$ exists; \ldots.

Manifolds of bounded geometry, but with empty boundary, have been
examined and used, among many others, by  the following authors:
Shubin \cite{Shubin(1992b)} studies (pseudo)differential operators, Roe
\cite{Roe(1988b)} derives an index theorem, Dodziuk \cite{Dodziuk(1981)} studies the Laplacian, Cheeger-Gromov-Taylor \cite{Cheeger-Gromov-Taylor(1982)} investigate integral kernels of functions of the Laplacian, Eichhorn deals with various of the above and other aspects \cite{Eichhorn(1989), Eichhorn(1991a), Eichhorn(1992a)}.

\begin{definition}\label{def_bounded} Suppose $M$ is a manifold with boundary $\boundary M$. It is of {\em bounded geometry} if the following holds:
\begin{itemize}
\item[(N)] {\em Normal collar:} there exists $r_C>0$ so that the geodesic collar 
\[ N=[0,r_C)\times \boundary M \to M: (t,x)\mapsto \exp_x(t\nu_x) \]
is a diffeomorphism onto its image, where $\nu_x$ is the unit inward normal vector at $x\in \boundary M$.
Equip $N$ with the induced Riemannian metric. Set $N_{1/3}:=\im[0,r_C/3)\times\boundary M$ and define {\em $N_{2/3}$} similarly.
\item [(IC)] {\em Positive injectivity radius of $\boundary M$:} $r_{inj}(\boundary M)>0$.
\item[(I)] {\em Injectivity radius of $M$:} There is $r_i>0$ so that for $x\in M-N_{1/3}$ the exponential map is a diffeomorphism on $B(0,r_i)\subset T_xM$.
\item[(B)] {\em Curvature bounds:} For every $K\in\N$ there is $C_K>0$ so that $\abs{\nabla^i R}\le C_K$ and $\abs{\bar\nabla ^i l}\le C_K$ $\forall 0\le i\le K$. Here $R$ is the curvature and $l$ second fundamental form tensor, $\nabla$ is the Levi-Civita connection of $M$ and $\bar\nabla$ the one of $\boundary M$
\end{itemize}
\end{definition}
This is a generalization of  corresponding definitions for manifolds
without boundary. The embedding of the boundary is described by the second fundamental form;  bounds on it will guarantee some homogeneity. Because the injectivity radius does not make sense near the boundary, we replace it by the geodesic collar.

Now we recall  Gaussian coordinates and their
substitute near the boundary.
\begin{definition}
Choose $0<r_i^C< r_{inj}(\boundary M)$ and $x'\in\boundary M$. Identify  $T_{x'}\boundary M$ with $\reals^{m-1}$ via an orthonormal base. Define {\em normal collar coordinates}
\[ \kappa_{x'}: \underbrace{B(0, r_i^C)}_{\subset \R^{m-1}}\times [0,r_C)\to M: (v,t)\mapsto \exp^M_{\exp_{x'}^{\boundary M}(v)}(t\nu) \]
(we compose the exponential maps of $\boundary M$ and of $M$).
Here $\nu$ is the inward unit normal vector field.\\
For $x\in M-N_{1/3}$ the exponential map yields {\em Gaussian coordinates}
\[ \kappa_x: B(0,r_i)\to M: v\mapsto \exp_x^M(v) \]
(identify $T_xM$ with $\R^m$ via an orthonormal base). We denote normal collar coordinates and Gaussian coordinates with the common name {\em normal coordinates} and the range $\kappa_x(B(0,\epsilon))$ (if $x\in M-N_{1/3}$) as well as $\kappa_{x}(B(0,\epsilon)\times [0,r_C))$ (if $x\in\boundary M$) with $NR(x,\epsilon)$.
\end{definition}

The main result of \cite{Schick(1998b)} is:
\begin{proposition} A Riemannian manifold $M$ has bounded geometry if
  and only if (N), (IC), (I) of \ref{def_bounded} hold and (instead of (B))
\begin{itemize}
\item[(B1)] For any $K\in\N$ an $M$-universal $C_K>0$ exists, so that in arbitrary normal coordinates (if $r_C$, $r_i$ and $r_i^C$ are sufficiently small) the following holds for the metric tensor $g_{ij}$ and its inverse $g^{ij}$:
\[ \abs{D^\alpha g_{ij}}\le C_K\quad\text{and}\quad \abs{D^\alpha g^{ij}}\le C_K\qquad\forall\abs{\alpha}\le K .\]
\end{itemize}
It suffices to check Condition (B1) for an atlas.
\end{proposition}

On our manifolds of bounded geometry we will study vector bundles, sections of these bundles and differential operators between them. All this has to be sufficiently regular:

\begin{definition} A bundle $E\bundleover M$ over a manifold of bounded geometry has bounded geometry if it is provided with the following additional structure:\\
a trivialization (called {\em admissible}) over each normal chart  so that the transition  functions for the trivializations, pulled back via normal charts, are bounded in the following sense: $\forall K\in\N$ an $E$-universal $C_K>0$ exists so that all partial derivatives of the transition functions up to order $K$ are bounded by $C_K$.
\end{definition}

\begin{notation} We will frequently talk of functions $f:M\to X$ ``in normal coordinates'' for  manifolds $M$ of bounded geometry. Then we always mean the composition of $f$ with the corresponding chart, but we omit writing down this chart. 
Similar remarks apply to sections of bundles of bounded geometry  and other objects like differential operators.
\end{notation}

\begin{definition}
Let $E\bundleover M$ be a bundle of bounded geometry. Suppose it is equipped with a Hermitian metric. We denote this metric {\em bounded} if the conjugate linear isomorphism 
$E\to E': x\mapsto (\cdot,x)$ is bounded with bounded inverse (in the sense of \ref{def_bounded_do}).

In an admissible trivialization the   Gram's matrix function
$b_{ij}(x)=(e_i,e_j)_x$
describes the Hermitian inner product. Let $b^{ij}(x)$ be the inverse. Boundedness translates to the existence of an  $E$-universal
 $C_K>0$  so that $|\partial_x^\alpha b_{ij}(x)|\le C_K$ and  $\abs{\partial_x^\alpha b^{ij}(x)}\le C_K$ $\forall \abs{\alpha}\le K$.
\end{definition}

\begin{example} (compare \cite{Shubin(1992b)} and \cite{Roe(1988b)})
\begin{enumerate}
\item\label{cpt} Every compact Riemannian manifold and every bundle
  over it are of bounded geometry.
\item If $\tilde M$ is the covering of a compact manifold $M$ and the
Riemannian metric is lifted then $\tilde M$ is of bounded
geometry. Every pullback bundle $\tilde E$ of a bundle over $M$ is of
bounded geometry (if we lift the trivializations).
\item Compact perturbations (e.g.~connected sums) do not change the bounded geometry property.
\item\label{fol} The leaves of a foliation of a compact Riemannian manifold with the induced metric are of bounded geometry.
\item\label{TM} On a manifold of bounded geometry, the trivial
as well as the tangent bundle (with trivializations given by normal coordinates) are  of bounded geometry.
\item If $M$ is of bounded geometry then the same is true for $\boundary M$
\item\label{comb} If $E,F$ are bundles of bounded geometry then so are $E\oplus F$, $E\otimes F$, $\Lambda^*E$, $\Hom(E,F)$, $E'=\Hom(E,\K)$, \ldots
\end{enumerate}
\end{example}

\begin{definition}\label{def_Cb}
 Let $E$ be a bundle of bounded geometry over $M$. Then
$C^\infty(E)$ denotes smooth sections, $C^\infty_0(E)$ smooth sections with compact support. Similarly defined are $C^k$ and $C^k_0$.
\[ C^k_b(E):=\{f\in C^k(E);\;|\partial^\alpha(t_p\circ f\circ \kappa_p)|\le C_k(f)\quad\forall\abs{\alpha}\le k\}.\]
Here, $\kappa_p$ runs through the normal charts, $t_p$ through the corresponding admissible trivializations.
$C^k_b(E)$ is a Banach space. The norm of $f$ is the smallest possible constant $C_k(f)$ in the definition.
\end{definition}

On non-compact manifolds, analysis is only possible if the operator classes and the function spaces match. Here, we will consider uniform function spaces and uniform operators. The prototype is the lift of an operator from a compact manifold to a covering.

\begin{definition}\label{def_bounded_do}
Suppose $E,F$ are bundles of bounded geometry over $M$. Let $A:C_0^\infty(E)\to C^\infty(F)$ be a differential operator of order $k$. In normal coordinates with admissible trivializations, $A$ becomes a matrix $(\sum_{|\alpha|\le k}a_\alpha^{ij}(x)D^\alpha)_{ij}$. We say $A$ is {\em bounded} if $\forall K\in\N$
\[ |D^\beta a^{i,j}_\alpha(x)|\le C_K\quad\forall|\beta|\le K\]
 with $(E,F)$-universal constants $C_K$. Similar notions apply to boundary differential operators $p:C_0^\infty(E)\to C^\infty(X)$, where $X$ is a bundle of bounded geometry over $\boundary M$. (A boundary differential operator is a differential operator composed with restriction to the boundary.)
\end{definition}

\begin{example}
On a manifold of bounded geometry the exterior differentiation $d$, its formal adjoint $\delta$, the Laplacian $\Delta$ and the Hodge-operator $*$ are bounded.
\end{example}
\begin{proof}
Boundedness is clear for
$d$. Since the composition of bounded operators is bounded and since
$\delta=\pm*d*$ and $\Delta=d\delta+\delta d$, we only have to check
it for $*$. But the entries in the matrix locally describing $*$ are
polynomials in the metric tensors $g_{ij}$ and $g^{ij}$.
\end{proof}

\forget{
\subsection{Technical properties of manifolds of bounded geometry}

\begin{lemma}\label{V}
Let $M$ be a connected non-compact $\boundary$-manifold of bounded geometry. We can find a function $V_0:\R_{\ge 0}\to \R_{\ge 0}$ with the following properties:
\begin{itemize}
\item $V_0$ is monotonous
\item $V_0(r)\to 0$ if $r\to 0$
\item $V_0(r)\to\infty$ if $r\to\infty$
\item $V_0(r)\le \vol(B(x,r))\qquad\forall x\in M$
\end{itemize}
\end{lemma}
\begin{proof}

\end{proof}
}

\subsection{Sobolev theory for manifolds of bounded geometry}

We will define Sobolev spaces using local data. For this, we need 
particular coverings and partitions of unity, whose existence is established in \cite[Proposition 3]{Schick(1998b)}:
\begin{lemma} \label{partition}
Let $M$ be a  manifold of bounded geometry. There is $r_m>0$ so that if $0<r<r_m$ then a countable covering of $M$ exists by charts $\{NR(x_i,r)\}_{i\in\Z}$ which has the following properties:
\begin{itemize}
\item $x_i\in\boundary M$ for $i\ge 0$, $x_i\in M-N_{2/3}$ for $i<0$.
\item There is $M_f<\infty$ so that $\forall s<r_m$ and $\forall x$ the
  intersection $NR(x,s)\cap NR(x_i,r)\ne\emptyset$ for at most $M_f$
  of the $x_i$.
\item $\{NR(x_i,r/2)\}_{i\in\Z}$ is a covering of $M$.
\end{itemize}
The charts are denoted $\kappa_i:B(0,r)\to NR(x_i,r)$ ($i<0$) and $\kappa_i:B(0,r)\times [0, r_C)\to NR(x_i,r)$ for $i\ge 0$.

To this covering, a subordinate smooth partition of unity $\{\varphi_i\}$ exists which is uniformly bounded, i.e.\ $\forall K\in\N$ we find $C_K>0$ so that 
\[ |D^\alpha\varphi_i|\le C_K\qquad\forall
i\in\integers\quad\forall\abs{\alpha}\le K\qquad\text{(in normal
coordinates).}\]
\end{lemma}

\begin{definition}\label{defHs}
Let $M$ be a bounded $\boundary$-manifold.
Choose data as in Lemma \ref{partition}. For $s\in\R$  define the {\em Sobolev norm}
\[
|f|^2_{H^s(M)}:=\sum_i|\varphi_if\circ\kappa_i|^2_{H^s(\kappa_i^{-1}(NR(x_i,r)))}\quad
\forall f\in C^\infty_0(M). \]
Here $|\cdot|_{H^s(\kappa_i^{-1}(NR(x_i,r)))}$ is the Sobolev norm on Euclidian space.

In the following, we will often abbreviate by writing $H^s(\R^m)$ or
$H^s$ instead of $H^s(\kappa_i^{-1}(\cdot))$ and $f$ instead of $f\kappa_i$.

Let $H^s(M)$ be the completion of $C^\infty_0(M)$ with respect to the Sobolev norm $|\cdot|_{H^s}$. In a similar way (using admissible trivializations) we define $H^s(E)$ if $E$ is a bundle of bounded geometry over $M$.

It is easy to prove that other choices yield equivalent Sobolev norms.
\end{definition}

The next proposition lists elementary properties of our Sobolev spaces, which follow from the corresponding properties of  $\reals^n_+$.
\begin{proposition}\label{Hs_prop}
Let $E$ be a bundle of bounded geometry over $M$. Suppose $X$ is a bounded vector bundle over $\boundary M$. Then the following holds for the Sobolev spaces $H^s(E)$ ($s,t\in \R$):
{\renewcommand{\labelenumi}{(\arabic{enumi})}
\begin{enumerate}
\item Sobolev embedding theorem: we have a bounded embedding
\[ H^s(E)\hookrightarrow C^k_b(E)\quad\text{whenever $s>m/2+k$} . \]
\item For $s<t$ we have a bounded embedding with dense image
\[ H^t\subset H^s. \]
This map is compact if and only if $M$ is compact. We define
\[H^\infty(E):=\bigcap H^s(E);\qquad H^{-\infty}(E):=\bigcup
H^s(E). \]
\item \label{Hcompare}
If $r,s,t\in\reals$ with $s<t$ and $\epsilon>0$ then there is
$C>0$ such that
\begin{equation*}
 \abs{\cdot}_{H^s}^2\le \epsilon \abs{\cdot}^2_{H^t}+C \abs{\cdot}^2_{H^r} .	       
	       \end{equation*}
\item\label{do_boundedness} Let $A:C^\infty(E)\to C^\infty(F)$ be a bounded differential operator of order $k$. Then $A$ extends to a bounded operator
\[ A: H^s(E) \to H^{s-k}(F)\quad \forall s .\]
\item Let $p:C^\infty(E)\to C^\infty(X)$ be a bounded boundary differential operator of order $k$. Then $p$ extends to an bounded operator
\[ p: H^s(E)\to H^{s-k-1/2}(X) \qquad\text{if } s>k+\frac{1}{2}.\]
\forget{
\item\label{Hermite_pairing} If $E$ is a bounded Hermitian or Riemannian bundle, then 
$H^s(E)$ and $H^{-s}(E)$ are dual to each other by extension of
\[ (f,g)=\int_M(f,g)_xd\,x\quad f,g\in C^\infty_0(E). \]
}
\end{enumerate}}
\end{proposition}


%% file: elliptic.tex
\section{Ellipticity and bounded geometry}\label{Ellipticity and bounded geometry}\label{sec_ell}

This section deals with elliptic boundary value problems. We define "uniform ellipticity" as the right condition for manifolds of bounded geometry. We obtain elliptic regularity results and an a priori estimate for solutions of uniformly elliptic boundary value problems, and we give a sufficient condition for essential self adjointness of elliptic  differential operators.

\forget{In this section, every manifold and bundle shall be bounded as well as every operator, if not stated otherwise.

Some of the computations to check uniformity are quite technical. However, if $M$ is a covering of a compact manifold and the operators are lifted, then all these results are  fulfilled trivially. 
} 

\begin{notation} If $M$ is a manifold of bounded geometry, $x'\in\boundary M$, then 
\[ Z(x',\epsilon):=K(\underbrace{B(x',\epsilon)}_{\subset\boundary M}\times[0,\epsilon))\]
is the cylinder with height and radius $\epsilon$ ($K$ is the geodesic
collar map). Note that $Z(x,\epsilon)\subset\R^m$ is mapped to $Z(x,\epsilon)\subset M$ under a normal collar chart.
\end{notation}

\begin{definition}
A bounded differential operator $A:C^\infty(E)\to C^\infty(F)$ of order $\mu\ge 1$ between $k$-dimensional bundles of bounded geometry is called {\em uniformly elliptic} if for its principal symbol $a_\mu(x,\xi)$ in admissible coordinates we have
\begin{enumerate}
\item the matrix $a_\mu(x,\xi)$ is invertible for all $(x,\xi)$ with $\xi\ne 0$ (ellipticity).
\item There is a constant $C$ so that 
$\norm{a_\mu^{-1}(x,\xi)}\cdot|\xi|^\mu\le C\quad\forall(x,\xi)$.
\end{enumerate}
We use the norms induced by admissible trivializations.
\end{definition}

Before we define elliptic boundary value problems we make a digression on conventions concerning arbitrary boundary value problems.

\subsection{Boundary value problems}

\begin{definition}\label{def_bvp}\label{def_comp_of_bvp}
 Let $E,F$ be vector bundles over a manifold $M$ and $X_0,\dots,X_r$ bundles over $\boundary M$. A system 
\[ (A,p_0,\dots,p_r): C_0^\infty(E) \to C^\infty(F)\oplus C^\infty(X_0)\dots\oplus C^\infty(X_r) \]
is called a {\em (differential) boundary value problem} (of order $\mu$) if $A$ is an ordinary differential operator (of order $\mu$) and $p_i$ are boundary differential operators. We set $\vect p=(p_0,\dots,p_r)$ and $\bvp P=(A,\vect p)$.\\
The composition of $(A,\vect p)$ and $(B,\vect q)$ (if this makes sense) is defined as $(A\circ B, (\vect q\circ A)\oplus \vect p)$.
\end{definition}

Next we go on to describe adjoints to boundary value problems. From now on, 
\[ \bvp P=(A,p_0,\dots,p_r): C^\infty(E)\to C^\infty(F)\oplus\bigoplus_\alpha C^\infty(X_\alpha)\]
 is a boundary value problem of order $\mu$ with order $p_\alpha=\alpha$ and
$\dim X_\alpha=n_\alpha$ (possibly $n_\alpha=0$), $\dim E=n=\dim F$ ($\alpha=0,\dots,r$).
\begin{definition}
We say $\vect{p}=(p_0,\dots,p_r)$ is {\em normal} if in all charts around $\boundary M$  the matrix $(a_{i,j}(x,\xi))_{ij}$ of the principal part of $p_\alpha$ has rank $n_\alpha$ for fixed $x$ and for $\xi=(0,\dots,0,\xi_m);\;\xi_m\ne 0$.

If $\vect{p}$ is normal and $n_\alpha=n\;\forall \alpha$, it is called a {\em Dirichlet system}.
\end{definition}
Normality says that the different boundary operators $p_i$ are independent of each other. Dirichlet systems are maximal normal systems. They have a unique adjoint and are unique up to transformations of lower order.

\begin{lemma}\label{Dirichlet_supplement}
 Suppose $\vect{p}=(p_0,\dots,p_r)$ is normal. Then a (normal) set of
boundary differential operators $\vect{t}=(t_0,\dots,t_r): E\to Y_0\oplus\dotsb\oplus Y_r$ exists so that $\vect p\oplus\vect t$ is a Dirichlet system.
\end{lemma}
\begin{proof}
Choose a Riemannian metric on $M$ and on $E$. Then
$p_i=A_i\partial^i_\nu + B_i + C_i$ where $\partial_\nu$ denotes
differentiation in inward normal direction, $B_i$ is a differential
operator of order $i$ of $\boundary M$ (a tangential differential
operator), $C_i$ is a boundary differential operator of order less
than $i$ and $A_i:E\to X_i$ is a linear bundle morphism which is
surjective by normality. Set $Y_i:=\ker A\subset E$ and
$t_i:=pr_{Y_i}\cdot\partial^i_\nu$ ($pr_{Y_i}$ being the orthogonal projection $E\to Y_i$).
\end{proof}

\begin{lemma}\label{Dirichlet_equivalence}
If $\vect p=(p_0,\dots,p_r)$ and $\vect q=(q_0,\dots,q_r)$ are two Dirichlet systems then $\Phi=(\Phi_{i,j})_{0\le i,j\le r}$ and $\Psi=(\Psi_{i,j})_{0\le i,j\le r}$ exist such that
\[ \vect p=\Psi\vect q;\quad \vect q=\Phi\vect p \]
(this implies $\Psi\circ\Phi=\1=\Phi\circ\Psi$) and such that
\[ \Phi_{i,j}:X_i\to Y_j\quad\text{and}\quad\Psi_{i,j}:Y_i\to X_j \]
are tangential differential operators  of order $j-i$ (i.e.\
$\Phi_{i,j}=0=\Psi_{i,j}$ for $i>j$) and so that
$\Phi_{i,i}$ and $\Psi_{i,i}$
are linear bundle isomorphisms.
\end{lemma}

\begin{proof}
First observe that it suffices to construct $\Phi$ and $\Psi$ locally (glue them together with a partition of unity). Next replace $\vect q$ by $\vect D: U\times \K^n\to U\times(\K^n)^{r+1}$ with $D_i=\1\cdot \partial_m^i$ (the general case follows by twice applying this).\\
It follows directly from the required properties of $\vect p$ that $\vect p=\Phi\vect D$ with appropriate $\Phi$. We construct $\Psi_{i,j}$ inductively in $j$:\\
For $j=0$ set $\Psi_{0,0}=p_0^{-1}$ ($p_0$ is invertible by definition of Dirichlet systems). For $i>0$ set $\Psi_{i,0}=0$.\\
Suppose $\Psi_{i,j}$ is already constructed for $j<j_0$. We have
\[ p_{j_0}= A_0\circ \partial_\nu^{j_0} + \sum_{i=1}^{j_0}A_i \partial_\nu^{j_0-i} \]
with $A_i$ tangential differential operators of order $\le i$ and $A_0$ invertible (by the Dirichlet property). Hence
\[ \begin{split} \partial_\nu^{j_0} &=A_0^{-1}p_{j_0} - \sum_{i=1}^{j_0}A_0^{-1}\circ A_i\circ \partial_\nu^{j_0-i}\\
&\stackrel{\text{ind.}}{=}A_0^{-1}p_{j_0} - \sum_{i=1}^{j_0}A_0^{-1}\circ A_i\circ\sum_{k=0}^{j_0-i}\Psi_{k,j_0-i} p_k .
\end{split} \]
Rearranging terms yields $\Psi_{k,j_0}$ with the required properties.
\end{proof}

\begin{lemma}\label{Dirichlet_solvability}
Suppose $\vect p$ is a Dirichlet system. For arbitrary sections
$\varphi_i\in C^\infty_0(X_i)$ we find $f\in C^\infty_0(E)$ with $p_i
f=\varphi_i$. 
\forget{Moreover, the support of $f$ may be confined in any open set containing the support of each of the $\varphi_i$.}
\end{lemma}
\begin{proof}
It suffices to consider the local case. Let $\vect p=\Phi\vect D$ and
$\vect D=\Psi\vect p$ as in Lemma \ref{Dirichlet_equivalence}. Let
$\vect\varphi=(\varphi_0,\dots,\varphi_{\mu-1})$\forget{ and
$\supp\vect\varphi=\bigcup\supp\varphi_i$. Set
$\vect\alpha=\Psi(\vect\varphi)$}. By Whitney's
\cite{Whitney(1934)} we find smooth $f$ with compact support\forget{
arbitrarily close to $\supp\vect \alpha$} so that $\vect D f=\alpha$. Then
\[ \vect p f=\Phi\vect D
f=\Phi(\vect\alpha)=\Phi\Psi(\vect\varphi)=\vect\varphi. \qquad\qed\]
\renewcommand{\qed}{}
\end{proof}
\forget{
\begin{remark}\label{Dirichlet_inf}
We can also define infinite Dirichlet systems
$(p_0,p_1,\dotsc)$. \ref{Dirichlet_equivalence} and
\ref{Dirichlet_solvability} remain true without change of the proofs.
\end{remark}
}

\begin{definition}\label{adjoint_bvp}
Suppose $(A,\vect p):C^\infty(E)\to C^\infty( F)\oplus\bigoplus_i
C^\infty(X_i)$ is a normal boundary value problem. Choose $\vect
t:C^\infty(E)\to\bigoplus_i C^\infty(Y_i)$ so that $\vect p\oplus\vect
t$ is Dirichlet\forget{ (\ref{Dirichlet_supplement})}. A similar
system
\[(B,\vect s\oplus \vect q):C^\infty(F)\to C^\infty(E)\oplus\bigoplus
C^\infty(X_i)\bigoplus C^\infty(Y_i)\]
 is called the {\em
formal adjoint} (with respect to given  metrics on $E$, $F$,
$X_i$, $Y_i$ and on $M$) if $\forall f\in C^\infty_0(E),\; g\in C^\infty_0(F)$
\begin{equation}\label{Green}
(Af,g)-(f,Bg) = (\vect p f,\vect s g)-(\vect t f,\vect q g).
\end{equation}
Here $(\vect pf,\vect sg):=\sum_i\int_{\boundary M}
(p_if(x),s_ig(x))_{X_i|_x} dx$, etc.
In particular $B$ is the formal adjoint of $A$.\\
We call $(B,\vect q)$ the {\em formal adjoint system} to
$(A,\vect p)$ (with respect to the{\em Greenian formula}
\eqref{Green}).
\end{definition}

If the manifold and the bundles are of bounded geometry and all
operators are bounded then by continuity the Greenian formula
\eqref{Green} extends to $f\in H^\mu(E)$, $g\in H^\mu(F)$.

\begin{theorem}
To every Dirichlet system $(A,\vect p,\vect t)$ a unique adjoint system $(B,\vect s,\vect q)$ exists.
\end{theorem}
\begin{proof}
Uniqueness: 
Only considering $f$ with $\supp f \cap \boundary M=\emptyset$ we see
that $B$ must be the unique formal adjoint of $A$. Let 
$(B,\vect s,\vect q)$ and $(B,\vect s_1,\vect q_1)$ be two adjoint systems. By
Lemma \ref{Dirichlet_solvability} and the Greenian formula:
\[ \begin{split} &(\varphi,\vect s g)-(\psi,\vect q g)=(\varphi,\vect s_1g)-(\psi,\vect q_1g)\quad\forall\varphi,\psi\in C^\infty_0\\
&\implies \vect s=\vect s_1\quad\text{and}\quad \vect q=\vect q_1.
\end{split}\]
Existence: 
First we construct the
adjoint locally. On a chart, integration by parts as in \cite[p.\
218]{Lions-Magenes(1972)}  yields
\[ (Af,g)_{L^2(F)}-(f,Bg)_{L^2(E)} =(\vect Df,\vect Ng)_{\boundary M} \]
for $g$ with support in a chart neighborhood. ($\vect D$ is defined
in the proof of \ref{Dirichlet_equivalence}). Now, by \ref{Dirichlet_equivalence} $\vect D=\Phi(\vect p\oplus\vect t)$ with a tangential  $\Phi$. Let $\Phi^*$ denote the formal adjoint of $\Phi$ (as an operator on $\boundary M$). Then
\[ (Af,g)-(f,Bg)=(\vect Df,\vect Ng)=(\vect p\oplus\vect
tf,\Phi^*\vect Ng) .\]
Let $\vect{s}$ and $\vect{q}$ be the corresponding components of
$\Phi^*\vect{N}$ which by uniqueness are globally defined. The
remaining conditions on them are checked as in \cite[p.\
218]{Lions-Magenes(1972)}.
\end{proof}


Now we will define uniformly elliptic boundary value problems. The idea is to copy the local definition (for elliptic boundary value problems) and to require uniformity so that local constructions patch together to global ones on manifolds of bounded geometry. 

\forget{
The definition of a (local) elliptic boundary value problem is the following (H\"ormander \cite[10.6.2]{Hoermander(1963)})
\begin{definition}\label{def_elbwp}
A boundary value problem $(A,p):C^\infty(E)\to C^\infty(F)\oplus C^\infty(X)$ is elliptic if the following conditions are fulfilled:
\begin{itemize}
\item $A$ is an elliptic operator of order $\mu$ 
\item $p$ is a normal system of boundary conditions
\item the homogeneous local frozen boundary value problem (transported to $\R^m$ via a chart)
\[  \sum A^0_{ij}(x_0,D)u_j=0\qquad\text{in }\R^m_{\ge 0}\qquad \sum p^0_{ij}(x_0,D)u_j=0 \qquad\text{in }\R^{m-1} \]
has no solution of the form
\[ u_j(x)=e^{i(x,\xi')}w_j(x_n) \qquad\text{with $\xi'\ne 0$}  \]
which is bounded in $\R^m_{\ge 0}$. ($A^0$ denotes the principal part). 
\item $\mu \dim E=2\dim X$ and the order of $p_{ij}$ with respect to $D_n$ is smaller than $\mu$.
\end{itemize}
\end{definition}

Unfortunately, it is not easy to see how one could require uniformity in this definition.
}
We use the existence of fundamental solutions for elliptic boundary value problems and pose uniformity conditions on these. In Proposition \ref{explicit_uniform} we will give a condition for uniform ellipticity which can be checked directly on the coefficients of the operators. 
\begin{definition}\label{def_elliptic_bvp}
Let 
\begin{equation*}(A,p_0,\ldots,p_r):C^\infty(E)\to C^\infty(F)\oplus
C^\infty(X_0)\oplus\cdots\oplus C^\infty(X_r)\end{equation*}
 be a differential boundary value problem of order $\mu\ge 1$, where $E,F$ are bundles of bounded geometry over $M$ and $X_i$ are bundles of bounded geometry over $\boundary M$.\\
We call $(A,\vect p)$ {\em uniformly elliptic} if it is an
elliptic boundary value problem in the sense of H\"ormander
\cite[10.6.2]{Hoermander(1963)}, if $A$ is  uniformly elliptic and if the following holds: Since $(A,\vect p)$ is elliptic, $\forall x\in\boundary M$ there is $0<\epsilon_x<r_C$ and a bounded fundamental solution \cite[10.4.1]{Hoermander(1963)}
\[ R_x:L^2(F|_{Z(x,\epsilon_x)})\oplus\bigoplus_{i=0}^r H^{\mu-i-1/2}(X_i|_{\boundary M\cap Z(x,\epsilon_x)}) \to H^\mu(E|_{Z(x,\epsilon_x)}) .\]
Set $ \bundlesum M^k:= H^k(F|_{Z(x,\epsilon_x)})\oplus\bigoplus_{i=0}^r H^{k+\mu-i-1/2}(X_i|_{\boundary M\cap Z(x,\epsilon_x)})$.\\
The fundamental solution fulfills
\begin{enumerate}
\item $(A,\vect p)\circ R_x(\vect F) \equiv\vect{F}$ on $Z(x,\epsilon_x)$
\item $R_x\circ (A,\vect p)u=u$ if $\supp u\subset Z(x,\epsilon_x)$  and $u\in H^\mu(E|_{Z(x,\epsilon_x)})$.
\end{enumerate}
For uniform ellipticity we require that it is possible to choose $\epsilon$ independently of $x$, and $R_x$ so that its norm is bounded uniformly in $x$. 
\end{definition}

\begin{remark}\label{ext_fund_def}
  In \cite[4.10, 4.11]{Schick(1996)} it is shown that for a uniformly
  elliptic boundary value problem not only the induced map $H^\mu\to
  L^2$ has uniform fundamental solutions but every map $H^{\mu+k}\to H^k$.
\end{remark}

\begin{example}
On a covering manifold $\tilde M\bundleover M$
the lift of any elliptic differential  boundary value problem is
uniformly elliptic.
\end{example}

We will now derive elliptic regularity in our context.
\begin{theorem} \label{elliptic_regularity}
Let $\bvp P=(A,\vect p): C^\infty(E)\to C^\infty(F)\oplus\bigoplus_i C^\infty(X_i)$ be uniformly elliptic of order $\mu$ (situation of \ref{def_elliptic_bvp}). If $u\in H^t(E)$  and $\bvp P(u)\in \bundlesum M^{s-\mu}$ (here $\bundlesum M^s=H^s(F)\oplus\bigoplus_i H^{s+\mu-i-1/2}(X_i)$) then
\[ u\in H^s(E)\quad\text{and }\abs{u}^2_{H^s(E)}\le C_{t,\bvp P} \cdot\left(\abs{\bvp{P}u}^2_{\bundlesum M^{s-\mu}} + \abs{u}^2_{H^t(E)}\right) \]
with $C_{t,\bvp P}$ independent of $u$.
\end{theorem}
\begin{proof}
By \cite[p.\ 270]{Hoermander(1963)} the hypotheses implies $u\in
H^s_{loc}(E)$, i.e.\ $\varphi u\in H^s$ for every smooth function
$\varphi$ with compact support. It remains only to check the norm
inequality. Take a normal covering (with centers $x_i$) with
subordinate admissible partition of unity $\{\varphi_i\}$ as in Lemma
\ref{partition} and so small, that we have local fundamental solutions
$R_i$. Let $M_f$ be the covering dimension.
  In the following we abbreviate $\abs{\cdot}_{H^s}$ and $\abs{\cdot}_{\bundlesum M^s}$ by $\abs{\cdot}_s$.
\begin{align*}
\abs{u}^2_s &= \sum_i\abs{\varphi_i u}^2_s =\sum_i\abs{R_i \bvp P\varphi_i u}^2_s \\
&\le C \sum_i\abs{\bvp P\varphi_i u}^2_{s-\mu} \qquad\text{(Remark \ref{ext_fund_def})}\\
&\le 2C\Bigl( \sum_i\abs{\varphi_i \bvp P u}^2_{s-\mu} +
\sum_i\bigl|\bvp P_i\underbrace{\sum_j\varphi_j}_{=1}
u\bigr|^2_{s-\mu}\Bigr)\qquad\text{($\bvp P_i:=\bvp
P\varphi_i-\varphi_i\bvp P$)}\\
&\le C'\Bigl(\abs{\bvp P u}^2_{s-\mu} + \sum_i\sum_j\abs{\bvp P_i\varphi_j u}^2_{s-\mu}\Bigr)\qquad\text{($C'=C 2^{M_f+2}$)}\\
\intertext{($\bvp P_i$ are boundary differential operators of order
strictly lower than $\bvp P$ with support in
$\supp\varphi_i$. Therefore, each sum $\bvp P_i\sum_j\varphi_j$ has
only $M_f+1$ non-trivial summands, and this, together
with $\abs{a+b}^2\le 2(\abs{a}^2+\abs{b}^2)$ yields the  factor $2^{M_f}$. As the derivatives of the coefficients of $\bvp P_i$ are bounded independently of $i$ we get:)}
&\le \underbrace{C'({M_f}+1)(1+\sup_i\norm{\bvp P_i}^2_{\boundedops(H^{s-1},\bundlesum M^{s-\mu})})}_{=:C_0}\Bigl(\abs{\bvp P u}^2_{s-\mu} + \sum_j\abs{\varphi_j u}^2_{s-1}\Bigr)\\
&\stackrel{\ref{Hs_prop} (\ref{Hcompare})}{\le} C_0\Bigl(\abs{\bvp
Pu}^2_{s-\mu} +
\frac{1}{2C_0}\underbrace{\sum_j\abs{\varphi_j u}^2_s}_{=\abs{u}^2_s}
+ C_t\underbrace{\sum_j\abs{\varphi_j u}^2_t}_{=\abs{u}^2_t}\Bigr) .
\end{align*}
The a priori estimate follows.
\end{proof}

The following proposition is proved
in \cite[Chapter 4]{Schick(1996)}.
\begin{proposition}\label{comp_of_elliptic}
The composition and the adjoint of  uniformly elliptic boundary value
problems is uniformly elliptic, as well.
\end{proposition}

Now we will derive regularity theorems for adjoint boundary value problems. We start with a theorem of H\"ormander \cite[Theorem 10.4.3 and p.\ 270]{Hoermander(1963)} which describes the local situation:
\begin{theorem}\label{Hoermander}
Let $M$ be a Riemannian manifold and $E,F\bundleover M$, $X_i\bundleover\boundary M$ be bundles with Riemannian metrics. Suppose $(A,\vect p)$ is an elliptic boundary value problem as in Definition \ref{def_elliptic_bvp}.\\
Suppose $f\in H^t_{loc}(F)$, $g\in H^s_{loc}(E)$ and $T_i\in H^{-\infty}_{loc}(X_i)$. If
\[ (f,A\varphi)+\sum_i(T_i,p_i\varphi) = (g,\varphi)\qquad\forall\varphi\in C^\infty_0(E) \]
then
$f\in H^{s+\mu}_{loc}$, $A^*f=g$ and $T_i\in H^{s+i+1/2}_{loc}$.\\
Here $A^*$ is the formal adjoint of $A$ and the pairings are extensions of the $L^2$-inner products for smooth sections with compact support.
\end{theorem}

We will need the following lemma to apply Theorem \ref{Hoermander}.
\begin{lemma}\label{adjoint_ref} 
Let $(A,\vect p)$ be an elliptic boundary value problem as in \ref{def_elliptic_bvp}. Suppose $f\in L^2(F)$, $g\in L^2(E)$ so that
\begin{equation}\label{orthog_eq}
 (Au,f)=(u,g)\quad\forall u\in C^\infty_0(E)\text{ with }\vect pu=0 .
\end{equation}
To arbitrary $u\in C^\infty_0(E)$ choose $u_i\in C_0^\infty(E)$ with
$p_i u_i=p_i u$ and $p_j u_i=0$ if $j\ne i$ (possible by
Lemma~\ref{Dirichlet_solvability} since $\vect p$ is normal by ellipticity). Set
\[ T_i(p_iu) = (u_i,g)-(Au_i,f). \]
Then $T_i\in H^{-\infty}_{loc}(X_i)$.
In particular, 
\begin{equation}\label{eq_allg}
 (Au,f)+\sum_{i=0}^r T_i(p_iu)= (u,g)
\end{equation}
and, by Theorem \ref{Hoermander}, $f\in H^\mu_{loc}$ and $A^*f=g$.
\end{lemma}
\begin{proof}
It follows directly from \eqref{orthog_eq} that $T_i$ is well defined. Obviously $T_i$ is linear. It remains to show that for an arbitrary real-valued $C_0^\infty$-function $\varphi$ (without loss of generality with support in some chart neighborhood with bundle trivialization) the functional $\varphi\cdot T_i$ is $H^N$-bounded for some $N\in\N$.

Extend $\vect p$ to a Dirichlet system $\vect p\oplus\vect s$. By
\ref{Dirichlet_equivalence} we find $\Psi$ and $\Phi$ so that
$\vect D=\Psi(\vect p\oplus\vect s)$ and $\vect p\oplus\vect
s=\Phi\vect D$. Set $\vect\alpha^i:=\Psi(0,\dots,0,\varphi
p_iu,0,\dots,0)$.  Via a chart we can work in Euclidian space. We
consider $\vect\alpha^i$ as a tuple of functions on $\R^{m-1}$.

Fix $v_k:[0,r_C)\to \R$ with $\partial^j_t v_k(t)|_{t=0}=\delta_{kj}$. Set
\[ \begin{split}
 &f_{\alpha^i}:= \sum_k\alpha_k^i(x')\cdot v_k(x_m);\quad x'\in\R^{m-1}\\
 &\implies D_kf_{\alpha^i}|_{x_m=0}=\partial_{x_m}^k f_{\alpha^i}|_{x_m=0} =\alpha^i_k
 \qquad \implies \vect D f_{\alpha^i}=\vect\alpha^i\\
 &\implies (\vect p\oplus\vect s)(f_{\alpha^i})=\Phi\vect
D(f_{\alpha^i})=(0,\dots,\varphi p_i u,\ldots).
\end{split} \]
Especially $f_{\alpha^i}$ can be used to compute $T_i(\varphi p_iu)$.
Now
\begin{align*}
\abs{T_i(\varphi p_iu)}^2&=\abs{(f_{\alpha^i},g)-(Af_{\alpha^i},f)}^2\\
&\stackrel{2}{\le} \abs{f_{\alpha^i}}^2_{L^2}\cdot\abs{g}^2_{L^2} + \abs{Af_{\alpha^i}}^2_{L^2}\cdot\abs{f}^2_{L^2}\\
& \stackrel{(\abs{g}^2+\abs{f}^2\norm{A}^2)}{\le}\abs{f_{\alpha^i}}^2_{H^\mu}
\stackrel{2^\mu}{\le}\sum_k\abs{\alpha^i_k(x')v_k(x_m)}^2_{H^\mu(\R^m)}\\
&= \sum_k\sum_{j=0}^\mu\int_0^{r_C}\abs{\partial_{x_m}^j\alpha^i_k(\cdot)v_k(x_m)}^2_{H^{\mu-j}(\R^{m-1})}dx_m\\
&\stackrel{(\sup_{j,k,x_m}\abs{\partial^j_{x_m}v_k(x_m)}^2)}{\le}\sum_{j,k}\int_0^{r_C}\abs{\alpha^i_k}^2_{H^{\mu-j}(\R^{m-1})}\\
&\stackrel{\mu r_C}{\le}\sum_k\abs{\alpha^i_k}^2_{H^\mu(\R^{m-1})} = \abs{\vect\alpha^i}^2_{H^\mu}\\
&= \abs{\Psi(0,\dots,\varphi p_iu,\dots,0)}^2_{H^\mu}\stackrel{\norm{\Psi}^2}{\le} \abs{\varphi p_iu}^2_{H^{2\mu}}\\
&\stackrel{\norm\varphi}{\le}\abs{p_i u}^2_{H^{2\mu}}.
\end{align*}
Note that the constants do not depend on $p_iu$. Therefore, we have
shown that $T_i\in H_{loc}^{2\mu}$.\\
Now \eqref{eq_allg} follows from \eqref{orthog_eq} and from $\vect p(u)=\vect p(\sum_i u_i)$.
\end{proof}

The regularity result allows the following  decomposition of $L^2$.
\begin{corollary}\label{L2_split}
Make the same assumptions as in Theorem \ref{Hoermander}.
Suppose $(A,\vect p)$ and $(B,\vect q)$ are elliptic and formally adjoint with respect to the Greenian formula
\[ (Ag,f)-(g,Bf)=(\vect p g,\vect sf)-(\vect t g,\vect qf)\qquad f\in C^\infty_0;\;g\in C^\infty_0.\]
Setting $\domain D_A^\infty:=\{g\in C^\infty_0;\;\vect p g=0\}$ we get an orthogonal decomposition
\begin{equation}\label{split1}
 L^2(F)=\underbrace{\{ f\in C^\infty\cap L^2;\;Bf=0,\;\vect q f=0\}}_{=:\ker B} \oplus \overline{A(\domain D_A^\infty)} .
\end{equation}
If the bundles are of bounded geometry and the boundary value problems are uniformly elliptic, the following stronger result holds: Let $\vect s^\infty=(s_\mu,s_{\mu+1},\dotsc)$ be a collection of boundary differential operators so that $\vect p^\infty:=\vect p\oplus \bigoplus \vect s^\infty$ is an infinite normal system. Then
\begin{equation}\label{split2}
 L^2(F)=\{ f\in H^\infty(F); Bf=0=\vect q f\} \oplus \overline{A\{g\in C^\infty_0(E);\vect p^\infty g=0\}} .
\end{equation}
\end{corollary}
\begin{proof}
The Greenian formula implies that $\ker B$ and $A(\domain D^\infty_A)$ in \eqref{split1} are orthogonal to each other.\\
Suppose  that $f\in L^2(F)$ is orthogonal to $A(\domain D^\infty_A)$. Then Lemma \ref{adjoint_ref} and Theorem \ref{Hoermander} show that $Bf=0$ and $f\in H^\infty_{loc}= C^\infty$. It remains to show that $\vect q f=0$:\\
Choose an arbitrary smooth boundary value $\vect t g$ with compact support. 
By Lemma~\ref{Dirichlet_solvability} we find $g\in C^\infty_0(E)$ with $\vect p g=0$ and $\vect t g$ as chosen. Let $\psi$ be a compactly supported cutoff function which is $1$ on $\supp g$. Then
\[\begin{split}
(\vect t g,\vect q f) & =(\vect t g,\vect q (\psi f)) =
(\underbrace{\vect p g}_{=0},\vect s(\psi f))+(Ag,\psi f)-(g,B(\psi
f))\\
 & = \underbrace{(Ag,f)}_{f\perp A\domain
D^\infty_A}-(g,\underbrace{Bf}_{=0}) =0 
\end{split}\] 
For the bounded geometry formula \eqref{split2} observe first $\ker B\subset H^\infty(E)$ because of the a priori estimates of Theorem \ref{elliptic_regularity}. Since $A:H^\mu(E)\to L^2(E)$ is bounded it suffices  to show that $\{g\in C^\infty_0(E);\;\vect p^\infty g=0\}$ is $H^\mu$-dense in $\{g\in C^\infty_0(E);\;\vect p g=0\}$. This follows from the following approximation result, proved in \cite[3.42]{Schick(1996)}:

If $E\bundleover M$ and $X_i\bundleover \boundary M$ are bundles of bounded geometry and 
 \[\vect p=(p_0,p_1,\dotsc):C^\infty(E)\to C^\infty(X_0)\oplus C^\infty(X_1)\oplus\dotsb\] is an infinite Dirichlet system then, to $\epsilon>0$ and 
\[ \vect u=(0,\dots,0,u_k,u_{k+1},\dotsc)\text{ with } u_i\in C^\infty_0(X_i)\text{ and }\supp u_i\subset K\]
($k\in\naturals$ and $K$ a fixed compact set), we find $u\in C^\infty_0(E)$ with 
\[ p_iu=u_i\quad \forall i\in\naturals\qquad \text{ so that } \abs{u}_{H^k(E)}<\epsilon. \qquad\qed \]
\renewcommand{\qed}{}
\end{proof}

On a manifold of bounded geometry, the above can be strengthened to a condition for (essential) self adjointness:
\begin{theorem}\label{L2_adjoint}
Let $(A,\vect p)$ be uniformly elliptic of order $\mu$ with uniformly elliptic formal adjoint $(B,\vect q)$ with respect to the Greenian formula
\[ (Af,g)-(f,Bg)=(\vect p f,\vect sg)-(\vect t f,\vect qg).\]
Set $\domain D_A:=\{f\in H^\mu(E);\;\vect pf=0\}$ and $\domain D_A^\infty:=\{f\in C_0^\infty(E);\;\vect p f=0\}$.\\
 Consider $A$ as an unbounded operator on $L^2(E)$ with domain $\domain D_A$. Denote the same operator but with domain $\domain D_A^\infty$ with $A_\infty$.\\
Then: The $L^2$-adjoint of $A$ and of $A_\infty$ is 
\[ A^t=A^t_\infty=B\text{ with domain }\domain D_B=\{g\in H^\mu(F);\;\vect qg=0\}.\]
Analogously, $A$ is the adjoint of $B$ and therefore the closure of $A_\infty$.
\end{theorem}
\begin{proof}
We get immediately that $\domain D_B$ lies in the domain of $A^t$ (and of $A_\infty^t$) and that $B$ and $A^t$ coincide on $\domain D_B$. Take therefore $f\in\domain D(A_\infty^t)\supset\domain D(A^t)$ with $g=A_\infty^t f$. We have to show: $f\in H^\mu(F)$, $\vect q f=0$ and $g=Bf$. By definition of the adjoint 
\[ (A\varphi,f)=(\varphi,g)\quad\forall\varphi\in C^\infty_0(E)\text{ with }\vect p\varphi=0.\]
Lemma \ref{adjoint_ref} now says $f\in H^\mu_{loc}(E)$ and $Bf=g$. Next we show $\vect q f=0$:
\[ \begin{split}
0 &= (A\varphi,f)-(\varphi,g)=(A\varphi,f)-(\varphi, Bf)\\
&\stackrel{\text{formal adjoint}}{=}(\vect p\varphi,\vect s f)-(\vect t\varphi,\vect q f) \stackrel{\vect p\varphi=0}{=} -(\vect t\varphi,\vect q f).
\end{split} \]
Note that by Lemma \ref{Dirichlet_solvability} we find $\varphi$ with $\vect p\varphi=0$ to arbitrary smooth $\vect t\varphi$ with compact support. Since $C^\infty_0$ is dense in $L^2$ it follows $\vect q f=0$.\\
Now we apply the a priori estimate of Theorem \ref{elliptic_regularity}:
\[ \abs{f}_{H^\mu}\stackrel{C}{\le} \abs{Bf}_{L^2}+\sum_i\abs{q_if}_{H^{\mu-i-1/2}}=\abs{Bf}_{L^2} <\infty .\]
This concludes the proof.
\end{proof}

We will see in the next section how Theorem \ref{L2_adjoint} shows
that the Laplacian with Dirichlet/Neumann boundary conditions is self adjoint.

Next, we have to recall further aspects of Hilbert space theory.

\begin{lemma} \label{Hilbert_space_isometries}
Let $A: \domain{D}\to H$ be a positive ($\ge c>0$) self adjoint unbounded
operator on a Hilbert space $H$. Then the natural Hilbert space norm
$|\cdot|_{\domain{D}}$ on ${\domain{D}}$ given by $|x|^2_{\domain{D}}=|x|_H^2+|Ax|_H^2$ is equivalent
to $|\cdot|_A$ with $|x|_A=|Ax|_H$. Moreover, $A: ({\domain{D}},|\cdot|_A)\to
(H,|\cdot|_H)$ is an isometry.
\end{lemma}

\begin{corollary}\label{equiv_norms}
Suppose $A: {\domain{D}}\subset L^2\to L^2$ is a strictly positive self adjoint differential operator of order $\mu$ and ${\domain{D}}\subset H^\mu$ is closed.
Then the norms $|\cdot|_A$ and $|\cdot|_{H^\mu}$ on ${\domain{D}}$ are
equivalent and 
 $A: {\domain{D}}\to L^2$ is an isometry if we equip ${\domain{D}}\subset H^\mu$ with $|\cdot|_A$.
\end{corollary}
\begin{proof}
By Lemma~\ref{Hilbert_space_isometries}, $({\domain{D}},|\cdot|_A)$ is
complete. Because $A$ is of order $\mu$,
$|\omega|_A=|A\omega|_{L^2}\stackrel{\norm{A}}{\le}
|\omega|_{H^\mu}\quad \forall\omega\in {\domain{D}}$. Since
$({\domain{D}},|\cdot|_{H^\mu})$ is  complete, too, the norms are equivalent by the open mapping theorem.
\end{proof}

\begin{theorem}\label{Hk_split}
Let $(A,\vect p)$ be a uniformly elliptic formally self adjoint
boundary value problem of order $\mu$.  Choose $\epsilon>0$  so that $(\epsilon+A,\vect p)$ is again uniformly elliptic
(this is possible by \cite[4.14]{Schick(1996)}). Define the boundary system
\[ \vect p^\infty =(\vect p, \vect p\circ (\epsilon +A), \vect p\circ (\epsilon+A)^2,\dotsc) \]
which is normal because $(\epsilon+A,\vect p)$ is elliptic. Set
\[ \domain D_{A^k} :=\{ f\in H^{k\mu}(E);\; \vect p f=0=\dots = \vect p\circ(\epsilon+A)^{k-1}f \}\quad k\in\naturals .\]
Suppose $A$ is a non-negative unbounded operator on $\domain D_A$. Then
\[ (\epsilon+A)^k:\domain D_{A^{k+l}} \to \domain D_{A^l}\]
 is an isomorphism of topological vector spaces and
\[ \domain D_{A^l} = \{f\in H^\infty(E);\; Af=0, \vect p f=0\} \oplus \overline{A\underbrace{\{f\in H^\infty(E);\;\vect p^\infty f=0\}}_{=:\domain D_\infty}}. \]
Moreover, $(\epsilon+A)^k$ respects the splitting.\\
The operator $(\epsilon+A)^k$ with domain $\domain D_{A^k}$ is the $k$-th power of $(\epsilon+A)$ with domain $\domain D_A$ in the sense of unbounded operators.
\end{theorem}
\begin{proof}
Observe $\domain D_{A^0}=L^2$. It suffices to prove the theorem for
$(\epsilon+A)^k:\domain D_{A^k}\to L^2$. Now $(\epsilon+A)^k$,
considered as an unbounded operator on $L^2$ with domain $\domain
D_{A^k}$, is self adjoint by Proposition \ref{comp_of_elliptic} and \ref{L2_adjoint}. Obviously, $(\epsilon+A)$ is positive on $\domain D_A$. Then so is its $k$-th power as unbounded operator, and $(\epsilon+A)^k$ with domain $\domain D_{A^k}$ is obviously a restriction of this operator. By elliptic regularity \ref{elliptic_regularity} they are equal. Moreover, $\domain D_{A^k}$ is the kernel in $H^{k\mu}$ of a bounded operator (namely a boundary differential operator of order $k\mu-1$) and hence closed in  $H^{k\mu}$. Therefore, Lemma \ref{equiv_norms} applies and yields the isomorphism. By Corollary \ref{L2_split} we have the desired splitting of $L^2(E)=\domain D_{A^0}$. Via the isomorphism $(\epsilon+A)^k$ we get corresponding splittings of $\domain D_{A^k}$. We  have to show only that the inverse images are actually as stated. First note that $(\epsilon+A)$ is an isomorphism of $\{f\in C^\infty(E)\cap L^2(E);\; Af=0, \vect p f=0\}$ onto itself. Consequently, the inverse image of this space is the space itself. For $\domain D_\infty$ note that the boundary conditions are designed in a way that $(\epsilon+A)^k$ maps $\domain D_\infty$ to itself. Suppose, on the other hand, that $f\in \domain D_{A^k}$ and $(\epsilon+A)^k f\in \domain D_\infty$. Then by elliptic regularity $f\in H^\infty(E)$ (note that $f\in D_{A^k}$ guarantees the required boundary regularity). The boundary conditions on $(\epsilon+A)^kf$ together with those on $f$ just say $f\in \domain D_\infty$.\\
We conclude that $(\epsilon+A)^k$ is an automorphism of $\domain D_\infty$ and then also of its closures.
\end{proof}

%% file: laplace.tex
\section{The Hodge decomposition}\label{sec_lap}

We start with a discussion of the proper boundary conditions for the Laplacian.

Let $M$ be a complete oriented Riemannian manifold with $\boundary
M=M_1\amalg M_2$ (possibly empty). Let $d$ be the differential on forms, $\delta=\pm*d$ its formal adjoint and $\Delta=d\delta+\delta d$ the Laplacian.

\begin{notation} For a  form $\omega$, $\omega|_{M_1}$ means the pullback to the submanifold $M_1$.
\end{notation}

\begin{definition}\label{def_Dirichlet}
 Set $\Omega^p(M)=C^\infty_0(\Lambda^pT^* M)$.

Define the boundary differential operators
\begin{eqnarray*} p_0: \Omega^*(M) \to \Omega^*(\boundary M)&:& \omega\mapsto (\omega|_{M_1}+ *\omega|_{M_2})\\
p_1: \Omega^*(M) \to \Omega^*(\boundary M)&:& \omega\mapsto  (\delta\omega|_{M_1}+ *d\omega|_{M_2}) \\
p_k: \Omega^*(M) \to \Omega^*(\boundary M)&:& \omega\mapsto 
\begin{cases} p_0( (1+\Delta)^{k/2}\omega) & \text{if }k\ge 2\text{ even}\\
p_1( (1+\Delta)^{(k-1)/2}\omega) & \text{if }k\ge 3\text{ odd}\end{cases}
\end{eqnarray*}

\[ \text{Set}\quad\vect p^k= p_0\oplus p_1 \oplus\dots\oplus p_{k};\quad k=0,1,2,\dotsc \]

Define now Dirichlet/Neumann boundary conditions ($k\ge 0$)
\begin{align*} bd(k+1)(\omega) &\quad:\iff\quad \vect p^j(\omega)\equiv 0\quad\forall j\le k ; \\
bd(\infty)(\omega) &\quad:\iff \quad bd(k)(\omega)\quad\forall k\in \N .
\end{align*}
\end{definition}

\forget{
\begin{remark}
For $\omega\in \Omega^p(M)$ we have
\[ \begin{aligned}
 bd(k)(\omega) \quad\iff\quad
\end{aligned}
\begin{aligned} \begin{cases}\omega|_{M_1}=0,\, \delta\omega|_{M_1}=0,\, \delta d\omega|_{M_1}=0,\ldots,\\
 *\omega|_{M_2}=0,\, *d\omega|_{M_2}=0,\ldots;
\end{cases}\end{aligned} \]
where we stop with exactly $k$ differentiations.
\end{remark}
}

\begin{notation}
Spaces of forms which fulfill boundary conditions are written in the following manner: 
\[ H^s(\Lambda^*(T^* M), bd(k)):=\{\omega\in H^s(\Lambda^*(T^* M));\,bd(k)(\omega)\} .\]
\end{notation}

The following statements are more or less well known. Proofs may be
found in \cite[Chapter 5]{Schick(1996)}.
\begin{proposition}\label{DGreen}
 On a complete oriented Riemannian manifold with boundary $\boundary M=M_1\amalg M_2$,
the boundary value problem $(\Delta,\vect p^1)$ is formally self adjoint in the sense of Definition \ref{adjoint_bvp} with respect to the Greenian formula
\[ (\Delta\omega,\eta)_{L^2(M)}- (\omega,\Delta\eta)_{L^2(M)} = (\vect p^1\omega,\vect t\eta)_{L^2(\boundary M)}-(\vect t\omega,\vect p^1\eta)_{L^2(\boundary M)} ,\]
where $\vect t$ is defined as follows 
($*$ denotes the Hodge-$*$-operator of $M$; $\hat*$ the one of $\boundary M$):
\[ \begin{split} t_0:\Omega^p(M)\to \Omega^*(\boundary M): &\omega\mapsto ((-1)^{pm}\hat{*}(*d\omega)|_{M_1} -\hat{*}(\delta\omega)|_{M_2})\\
 t_1:\Omega^p(M)\to \Omega^*(\boundary M): &\omega\mapsto ((-1)^{m(p+1)}\hat{*}(*\omega)|_{M_1}-\hat{*}\omega|_{M_2}) .
 \end{split} \]
\end{proposition}

\begin{definition}\label{def_of_harmonic}
Define the space $\mathcal{H}^p( M, M_1)$ of harmonic forms which fulfill {\em absolute boundary conditions} on $ M_1$ and {\em relative boundary conditions} on $ M_2$ as
\[ \mathcal H^p( M, M_1):= \{\omega\in C^\infty\cap L^2(\Lambda^pT^* M);\, \Delta\omega=0,\,bd(2)(\omega)\} .\]

If $M$ is complete then these forms are closed and coclosed, i.e.\
$\mathcal H^p(M,M_1)= \{\omega\in C^\infty\cap L^2(\Lambda^pM);\,d\omega=0=\delta\omega, \omega|_{ M_1}=0= *\omega|_{ M_2}
\}$.
\end{definition}

As application of Corollary \ref{L2_split} we obtain:
\begin{theorem}\label{L2Hodge}
Let $M$ be a complete oriented Riemannian $\boundary$-manifold. Then we have the orthogonal {\em Hodge decomposition}:
\[ L^2(\Lambda^p T^*M)= {\mathcal H^p}(M, M_1)\oplus \overline{d\Omega^{p-1}_d} \oplus\overline{\delta\Omega^{p+1}_\delta} , \]
where $\Omega_d^p:=\{\omega\in C^\infty_0;\;\omega|_{M_1}=0\}$ and $\Omega_\delta^p:=\{\omega\in C^\infty_0;\;(*\omega)|_{M_2}=0\}$.
\end{theorem}

\subsection{Hodge decomposition for manifolds of bounded geometry}
 We want  to show that on a manifold of bounded geometry the Dirichlet/Neumann boundary value problem for the Laplacian $(\Delta,\vect p^1)$ is uniformly elliptic. 

To do this we investigate how to check for uniform ellipticity directly from the coefficients of the boundary value problem, without explicitly constructing the fundamental solutions.
 
\begin{definition}\label{expl_uni}
Let $\bvp P=(A,\vect p)$ be a boundary value problem as in Definition \ref{def_bvp}. For each point $b\in\boundary M$ we get normal collar coordinates and admissible bundle trivializations. In these  coordinates the problem is described by matrices $(\sum_\alpha A^{(ij)}_\alpha(x)D^\alpha)_{ij}$  and $(\sum_\alpha p^{(ij)}_{\alpha,k} D^\alpha)_{ij}$ $(k=0,\dots,\mu-1)$ with corresponding symbols $A(x,\xi),\dotsc$ (omiting an additional index $b$). In the following we inspect the exposition in H\"ormander's \cite[chapter 10]{Hoermander(1963)}. We consider only the principal parts $A^0$ and $p_k^0$. Substituting $D=\partial/\partial_t$ for $\xi_m$ in the principal symbol and freezing $x$ and $\xi'=(\xi_1,\dots,\xi_{m-1})$ we get an ordinary boundary value problem on $\R_{\ge 0}$. As our boundary value problem is elliptic we find $C_{x,\xi}>0$ so that $\forall v=(v_1,\dots,v_n)\in C^\infty(\R,\K^n)$ (compare \cite[10.2.3]{Hoermander(1963)}):
\begin{multline*} \sum_l\sum_{i=0}^m\int_0^\infty\abs{D^iv_l}^2dt + \sum_l\sum_{j=0}^{m-1}\abs{D^jv_l(0)}^2\\
 \le C_{x,\xi'}\Bigl( \sum_{l'}\int_0^\infty\abs{( A^0(x,\xi',
D)v)_{l'}}^2 dt +\sum_{k,l'}\abs{(p_k^0(x,\xi',D)v)_{l'}(0)}^2\Bigr) .
\end{multline*}
The smallest possible constant depends continuously on the  matrices. Let $C^b_x<\infty$ denote the best possible constant valid $\forall\abs{\xi'}=1$.\\
Moreover, if $\lambda_{s,\xi'}-1$ is a bound for the absolute value of the zeros of the polynomial $p(z)=A^0(x,\xi',z)$ then $C^{\lambda_{x,\xi'}}_{x,\xi'}>0$ exists so that
\begin{multline*}
 \sum_l\sum_{i=0}^m\int_0^\infty e^{2\lambda t} \abs{D^jv_l}^2\,dt+\sum_l\sum_{j=0}^{m-1}\abs{D^jv_l(0)}^2\\
 \le C^\lambda_{x,\xi'}\sum_{l'}\int_0^\infty e^{2\lambda t}\abs{(A^0(x,\xi',D)v)_{l'}}^2\,dt  .
\end{multline*}
Again, $C^\lambda_{x,\xi'}$ depends continuously on $\lambda$ and the matrices. For $\lambda\ge\lambda^b_x:=\max\{\lambda_{x,\xi'};\;\abs{\xi'}=1\}$  set $C^b_{x,\lambda}:=\max_{\abs{\xi'}=1}\{C^\lambda_{x,\xi'}\}$.

Let $\{c_i^b(x)\}_i$ be the collection of all coefficients of the principal part of any of the differential operators in question (this is for fixed $b$ only a finite set). Set
\[ \delta_b=\frac{1}{\sum_i\sup_{\abs{x}\le r}\abs{c_i(x)}} .\]
Let $r$ be a Euclidian radius so that the via a normal chart transported boundary value problem is defined on $B(0,r)\subset\R^m$.
\end{definition}

\begin{proposition}\label{explicit_uniform}
If $\lambda,r,C>0$  exist so that (notation as above)
\begin{itemize}
\item[(C0)] $\lambda^b_x\le \lambda\quad\forall x,b$
\item[(C1)] $C_{x,\lambda}^b\le C$ and $C_x^b\le C\quad\forall x,b$
\item[(C2)] $\abs{c_i^b(x)}\le C\quad\forall i,b,\abs{x}\le r$
\end{itemize}
then we can find $\rho$ independently of $b$ such that in normal coordinates around $b$ a fundamental solution $R_b$ on $B(b,\rho)$ in the sense of Definition \ref{def_elliptic_bvp} exists and the norm of $R_b:\bundlesum M_0\to H^\mu(\K^n\bundleover\R^m_+)$ is bounded by a constant independent of $b$.
\end{proposition}
\begin{proof}
Close inspection of H\"ormander's construction \cite[Theorem 10.4.1]{Hoermander(1963)} of the fundamental solution shows that for the domain of the fundamental solution we take a ball with radius $\rho_b=C_b\cdot \delta_b$ where $C_b$ depends only on $C_x^b$, $\lambda$ and $C^b_{x,\lambda}$. Following H\"ormander's construction further we see that the norm of this fundamental solution is bounded by a universal expression in $C_x^b$, $\lambda$ and $C^b_{x,\lambda}$. 
Considering the definition of $\delta_b$ this finishes the proof. 
\end{proof}

\begin{proposition}\label{unif_Laplace}
On a manifold of bounded geometry the boundary value problems $(\Delta,\vect p^1)$ and $(1+\Delta, \vect p^1)$ are uniformly elliptic.
\end{proposition}
\begin{proof}
It is well known that these boundary value problems are elliptic (compare Schwarz \cite{Schwarz(1995)}). The principal symbol of $\Delta$ on forms in normal coordinates around $b\in M$ is just $(\sum_{i,j}g_b^{ij}(x)\xi_i\xi_j)\cdot\1$. To prove uniform ellipticity of the Laplacian we must show that 
\[ \frac{\abs{\sum g_b^{ij}(x)\xi_i\xi_j}}{\sum_i\abs{\xi_i}^2};\quad\xi\ne 0 \]
is bounded from below by some positive constant. Since the matrix $g=(g^{ij})_{ij}$ is symmetric and positive and by substituting  $\sqrt{g}^{-1}\xi$ for $\xi$ we see that this is equivalent to finding an upper bound for
\[ \frac{\abs{\sum g_{ij}^b(x)\xi_i\xi_j}}{\sum\abs{\xi_i}^2} \]
where $g_{ij}$ is the inverse of $g^{ij}$. By the  definition of bounded geometry such a bound exists independently of $x$ and $b$.

It remains to check boundary uniformity. Here we will use Proposition \ref{explicit_uniform}. What is the local expression for $\Delta$ and $\vect p^1$? We will not write down the complicated formulas (see Schwarz \cite{Schwarz(1995)}) but observe that the coefficients $c_i^b$ are polynomials in $g_{ij}^b$, $g^{ij}_b$ and their derivatives up to order 2. Therefore, condition (C2) of Lemma \ref{explicit_uniform} is fulfilled. Concerning (C0) and (C1) we have the map 
\[\begin{CD} M\times \overline{B(0,r)}\times S^{m-2} @>>> \K^R @>>> \R\times\R\times\R \\
 (b,x,\xi') @>>> (c^b_i(x,\xi'))_i @>>> (C_x^b,\lambda,C^b_{x,\lambda})
\end{CD}\]
where in our specific case each of the $c_i^b(x,\xi')$ is a polynomial in $\xi'$ and derivatives of $g_{ij}^b(x)$ and $g^{ij}_b(x)$. The map $\K^R\to \R^3$ is only defined on the subset $E\subset \K^R$ which comes from elliptic boundary value problems and is continuous on this subset. We will show that the closure of the range of $M\times B(0,r)\times S^{m-2}$ is a compact subset of $E$. Then on this set
bounds  for $\lambda^b_x$ and for 
$C_x^b$ and $C_{x,\lambda}^b$ exist and the proof is finished.\\
For compactness, take a sequence $(b_n,x_n,\xi'_n)\subset M\times B\times S^{m-2}$. We have to produce a subsequence so that the image converges to some element of $E$. But we may choose a subsequence so that (after relabeling) $(\xi'_n)$, $x_n$, $g_{ij}^{b_n}(\cdot)$ and $g^{ij}_{b_n}(\cdot)$ converge for each $i,j$ (the latter in $C^2$-norm on $\overline{B(0,r)}$ by Ascoli's theorem and uniform boundedness (of the derivatives up to order $3$ of $g^{**}$ and $g_{**}$)) to $\xi'$, $x$ and $g_{ij}$, $g^{ij}$ respectively. Each of the algebraic relations between the $g_{ij}^n$ and $g^{ij}_n$ are preserved under the limit. Especially, $g_{ij}$ and $g^{ij}$ form the components of the covariant and contravariant metric tensor of some twice differentiable Riemannian metric $g$ on $B(0,r)$. Similarly, $(c_i(\xi'))_i:=\lim_n c_i^{b_n}(x_n,\xi'_n)$ are just the corresponding coefficients at $x$ for the boundary value problem $(\Delta,\vect p^1)$. In particular, they come from an elliptic boundary value problem, i.e.\ they lie in $E$.
\end{proof}

Now we can apply our theory of uniformly elliptic (and formally self adjoint) boundary value problems to $(\Delta,\vect p^1)$. It only remains to check:
\begin{lemma}\label{Dpos1}
The operator $\Delta$ is non-negative on its domain $\domain D_\Delta=\{\omega\in H^2(\Lambda^pT^*M);\;\vect p^1\omega=0\}$. More precisely:
\begin{equation}\label{positivity} (\Delta\omega,\omega)=(d\omega,d\omega)+(\delta\omega,\delta\omega)\ge 0 \qquad\forall\omega\in\domain D_\Delta.
\end{equation}
\end{lemma}
\begin{proof}
Simply integrate by parts. By continuity, this works for $\omega\in H^2$ because $d$ and $\delta$ are bounded operators in the sense of bounded geometry and therefore bounded from $H^2$ to $H^1$ to $L^2$.
\end{proof}

Theorem \ref{elliptic_regularity} implis $\mathcal H^p(M,M_1)\subset H^\infty(\Lambda^p T^*M)$. 
Applying  \ref{Hk_split} to our situation, we get:
\begin{theorem}\label{HsHodge}
Let $M$ be an orientable manifold of bounded geometry, possibly with boundary. Then
\begin{equation*}
\domain D_{\Delta^k}:=H^{2k}(\Lambda^p(T^*M), bd(2k))={\mathcal H^p}(M, M_1) \oplus \overline{d\Omega^{p-1}_\infty} \oplus \overline{\delta\Omega^{p+1}_\infty} .
\end{equation*}
Here  $\Omega_{\infty}^p:=H^\infty(\Lambda^p(T^*M),bd(\infty))$.
The closure is taken with respect to the given topology on $H^{2k}$, and the decomposition is orthogonal with respect to the Hilbert space structure pulled back from $L^2$ via the following isometry (it induces the $H^{2k}$-topology): 

The unbounded operator $(1+\Delta)^k$ on $L^2(\Lambda^*(T^*M))$ with domain $\domain D_{\Delta^k}$ is positive self adjoint and considered as an operator 
\[ (1+\Delta)^k: \domain D_{\Delta^k} \to L^2(\Lambda^*(T^*M)), \]
is an isometry which respects the decomposition above.

$(1+\Delta)^k$ with domain $\domain D_{\Delta^k}$ is the $k$-th power in the sense of unbounded operators of $(1+\Delta)$ with domain $\domain D_\Delta$.
\end{theorem}


%% file: L2+curvature.tex
\section{$L^2$-cohomology and curvature}\label{sec_curv}

In this section we examine the relations between $L^2$-de Rham cohomology and curvature. This is based on the Weizenb\"ock formula $\Delta=\nabla^*\nabla-\mathcal R^W$ with its integrated consequence (for $\boundary$-manifolds) 
\begin{equation}\label{eq_Bochner}\begin{split} (\nabla\omega,\nabla\eta)_{L^2}&=(\mathcal R^W\omega,\eta)+(d\omega,d\eta)+(\delta\omega,\delta\eta)+\\
&\quad\int_{M_1}(\mathcal
S\omega,\eta)_{\Lambda^*(T_x^*M)}\,dx+\int_{M_2}(\mathcal S{*\omega},*\eta)_{\Lambda^*(T_x^*M)}\,dx
\end{split}\end{equation}
for $\omega,\eta\in C_0^\infty(\Lambda^pT^*M)$ with $\omega|_{M_1}=0=\eta|_{M_1}$ and $*\omega|_{M_2}=0=*\eta|_{M_2}$ (proved in Schwarz
\cite[2.1.5 and 2.1.7]{Schwarz(1995)}).
Here the (symmetric) endomorphism $\mathcal S_k\in\End(\Lambda^k(M)|_{\boundary M})$ is defined by the identity 
\begin{equation}\label{Weiz}
(\mathcal
S\omega,\omega)_{\Lambda^k}=-(\omega,i_\nu d\omega)_{\Lambda^k}+(\delta\omega,i_\nu\omega)_{\Lambda^{k-1}}-1/2
\frac{\partial}{\partial\nu}(\omega,\omega)_{\Lambda^k},
\end{equation}
 where $\nu$ is the unit inward normal field, and $i_X$ denotes contraction of a form with the
vector field $X$.

We derive algebraic properties of the operators $\mathcal S$ and $\mathcal R^W$, especially sufficient conditions for negativity. With this at hand, it is an easy task to get vanishing results for square integrable harmonic forms.\\

\begin{proposition}\label{pos_def}
Let $M$ be a Riemannian manifold, $R^W_p$ the Weizen\-b\"ock tensor and $S_p$ the boundary tensor \ref{Weiz}.
\begin{enumerate}
\item\label{R1}\label{Rsec}
$\mathcal R^W_1=-\Ric$ ($\Ric$ the Ricci operator). If the curvature operator is $\ge 0$ then $\mathcal R^W_p\le 0$ $\forall p$.
\item\label{Sp}
$\mathcal S_p|_x\le 0$ if the sum of the $m-p$ smallest eigenvalues of $l|_x$ is non-negative ($l$ the second fundamental form operator with respect to the inward unit normal). Especially, if $\trace\,l\ge 0$, then $\mathcal S_1\le 0$; and if the second fundamental form is $\ge 0$ then $\mathcal S_p\le 0$ $\forall p$.
\end{enumerate}
\end{proposition}
Note: Negativity of the second fundamental form measures convexity of
the boundary inside
the manifold.

\begin{proof}[Proof of Proposition \ref{pos_def}]
(\ref{R1}) is classical (compare \cite[II.8.3 and II.8.6]{Lawson-Michelsohn(1989)}). For (\ref{Sp}) we explicitly perform the computations:\\
We consider the second fundamental form as a symmetric fiber-wise
operator $l: T\boundary M\to T\boundary M$. Let $e_1,\ldots,e_{m-1}$ be an orthonormal base for $T_b\boundary M$ of eigenvectors of $l$ with
eigenvalues $\lambda_i$. Then $\nu,e_1,\dots,e_{m-1}$ is an
orthonormal base for $T_bM$ at $b\in\boundary M$. Let $\omega$ be the $p$-form 
\[\omega_b=\sum_{i_1<\dots<
i_p}\alpha_{i_1\dots i_p}e_{i_1}\wedge\dots\wedge e_{i_p} + \sum\beta_{i_1\dots i_{p-1}}\nu\wedge
e_{i_1}\wedge\dots\wedge e_{i_{p-1}}.\]
 Then the explicit formulas of Schwarz \cite[(2.1.13)]{Schwarz(1995)} say that the coefficient in $\mathcal S\omega_b$ of $e_{i_1}\wedge\dots\wedge
e_{i_p}$ is $0$ and the coefficient $\gamma_{i_1\dots i_{p-1}}$ of $\nu\wedge e_{i_1}\wedge\dots\wedge e_{i_{p-1}}$ is
\begin{equation*}\begin{split}
&\gamma_{i_1\dots i_{p-1}} = (\mathcal S\omega)_b(\nu,e_{i_1},\dots,e_{i_{p-1}})=-\sum_r(l^\Lambda(e_r)\omega)(e_r,e_{i_1},\dots,e_{i_{p-1}}) 
\\
&\text{with}\quad(l^\Lambda(e_r)\omega)(e_r,e_{i_1},\dots,e_{i_{p-1}}) \\
&\qquad := \omega( \underbrace{l(e_r,e_r)}_{\lambda_r\nu},e_{i_1},\dots,e_{i_{p-1}}) + \sum_{l=1}^{p-1}\omega(e_r,\dots,\underbrace{l(e_r,e_{i_l})}_{0\text{ if $r\ne i_l$}},\dotsc)\\
&\qquad = \begin{cases}0&;\text{ if } r\in\{i_1,\dots,i_{p-1}\}\\
  \lambda_r\beta_{i_1,\dots,i_{p-1}}&; \text{ else}.
  \end{cases}
\end{split}\end{equation*}
Therefore, $
\gamma_{i_1,\dots,i_{p-1}} =-(\sum_{r\notin\{i_1,\dots,i_{p-1}\}}\lambda_r)\beta_{i_1,\dots,i_{p-1}}$. 

Hence $(\mathcal
S_p\omega,\omega)_b=-\sum_{i_1<\dots<i_{p-1}}(\sum_{r\notin\{i_1,\dots,i_{p-1}\}})\abs{\beta_{i_1,\dots,i_{p-1}}}^2$,
and $S_p|_b\le 0$ if and only if the sum of the smallest $(m-1)-(p-1)=m-p$ of the eigenvalues $\{\lambda_r\}$ of $l|_b$ is non-negative.
For $p=1$, this is the trace.
\end{proof}

\begin{theorem} ($L^2$-Bochner)\label{Boch}\\
Let $M^m$ be a complete orientable manifold with boundary. Suppose $M$ has infinite volume.
\begin{itemize}
\item If the Weizenb\"ock endomorphism $\mathcal R^W$ and the fundamental form  $\mathcal S$ are both
negative semidefinite on $k$-forms, then $\mathcal{H}^k(M,\boundary M)$ vanishes. If $\mathcal S_{m-k}\le 0$ instead of $\mathcal S_k$, then $\mathcal{H}^k(M)=\{0\}$.
\item
If the Ricci tensor of $M$ and the trace of the second fundamental form $l$ of $\boundary
M\subset M$ are both $\ge 0$, then $\mathcal{H}^1( M, \boundary M)$ vanishes.
\end{itemize}
\end{theorem}

\begin{proof}
First note that the second statement is a direct consequence of the first by Proposition \ref{pos_def}.
The first statement is proved  similar to the way 
Dodziuk treats manifolds without boundary in \cite{Dodziuk(1981b)}, using equation \eqref{eq_Bochner} (details are given in \cite[Theorem 7.7]{Schick(1996)}). 
\end{proof}

One can apply this to infinite coverings
of a compact manifolds. Using the $L^2$-Hodge de Rham Theorem \ref{Hodge_de_Rham_theorem}, the conclusion is:\\
If the base manifold carries a metric which fulfills the positivity conditions in the
$L^2$-Bochner theorem, then the corresponding $L^2$-cohomology groups $H^*_{(2)}(M,M_1)$ vanish.

\begin{example}
Take $\R P^3\#\R P^3\#\R P^3$ or $(\R P^3\#\R P^3\#\R P^3)-D^3$. An explicit calculation of Lott
and L\"uck \cite[theorem 0.1]{Lott-Lueck(1995)} shows that the first $L^2$-cohomology of both of them is not zero. Therefore they can not
carry a metric with non-negative Ricci curvature and (in the presence of boundary) non-negative trace of the second
fundamental form. (This does not follow from the classical  Bochner theorem.)
\end{example}

In the case of normal coverings of compact manifolds (and for other
manifolds of bounded geometry with sufficient symmetry), we can give
much weaker condition for the vanishing of $L^2$-cohomology
groups. On  a non-compact manifold the topology
does not restrict the geometry very much. For example, $\R^m$ can be
equipped with the flat Euclidian metric or with the hyperbolic metric
with constant negative curvature, which have quite different (rough)
qualities.

Therefore, it is natural to restrict the metrics to consider. We will
look for metrics in the bilipschitz class of a given metric only. This
is proposed in Roe \cite{Roe(1988c)}. Recall the following:

\begin{definition} Let $(M,g)$ and $(N,h)$ be two Riemannian manifolds, $f:M\to N$ a diffeomorphism. $f$ is called {\em bilipschitz}, if $C>0$ exists so that for the norms the following holds: $\norm{T_xf}\le C\;\forall x\in M$ and $\norm{T_yf^{-1}}\le C\;\forall y\in N$.
\end{definition}
Note that all the metrics on a compact manifold are bilipschitz. The same is true for the lifts to any covering.

We will strengthen the $L^2$-Bochner theorem by imposing positivity
conditions on the Weitzenb\"ock endomorphism only on (the 
small) eventually large sets (Definition \ref{eventlarge}).
Observe that 
  in a Riemannian manifold with infinite diameter, the complement of
  any compact set is eventually large.

\begin{theorem}\label{str_Boch}
Let $(M,g_1)$ be an oriented manifold of bounded geometry (possibly
with boundary), so that $\mathcal R^W_p\le 0$ on an eventually large
set and additionally $\mathcal S_p\le 0$ on an eventually large set.
Let $g$ be a metric on $M$ which is bilipschitz to $g_1$ and so that
$(M,g)$   admits a cocompact group $\Gamma$ of isometries.

Then $\mathcal H^p_{(2)}(M,\boundary M)=\{0\}$. If the condition
$\mathcal S_p\le 0$ is replaced by the condition $\mathcal S_{m-p}\le
0$, then $\mathcal H^p_{(2)}(M)=\{0\}$.
\end{theorem}
\begin{remark}
If $(M,g)$ is an infinite normal covering of a compact manifold it
obviously fulfills the conditions of the theorem.  
\end{remark}

The new and surprising feature  is the existence of an enormous subset where
the curvature assumption is not required to hold. This is due to the
uniformity of the metric  and
does not hold in general for manifolds of bounded geometry as shows \cite[Example 7.2]{Schick(1996)}.

For the proof of Theorem \ref{str_Boch} we need the following proposition which is essentially due to Roe \cite[1.11]{Roe(1988c)} (if $\boundary M=\emptyset$). We denote with $\Delta$ the unbounded operator on $L^2$ with domain $\domain D_\Delta:=\{\omega\in H^2(\Lambda^* T^*M);\;\omega|_{\boundary M}=0=\delta\omega|_{\boundary M}\}$ (compare \ref{def_Dirichlet}). Let $P$ be the orthogonal projection onto its kernel $\mathcal H^p(M,\boundary M)$.
\begin{proposition}\label{roe}
Suppose $M$ is a $\boundary$-manifold of bounded geometry. If
$\mathcal R^W_p$ and $\mathcal S_p$ are non-positive on an eventually
large subset $X$ of $M$, then $\forall \epsilon>0$ an eventually large subset $Y\subset M$ exists so that
\[ \sup_{x\in Y}\{\abs{Ps(x)}\}\le \epsilon\abs{s}_{L^2}\quad\forall s\in L^2(\Lambda^p T^*M) .\]
\end{proposition}
\begin{proof}
Set $K:=M-X$. 
Suppose $s\in L^2$ and set $s_t:=e^{-t\Delta}s$. Then
$s_t\stackrel{t\to\infty}{\to} Ps$ in $L^2$ by the spectral theorem
and $s_t\in\domain D_\Delta$. Even better: $s_t\in \domain
D_{\Delta^k}$, the domain of the $k$-th power of $(1+\Delta)$,
$\forall k\in\N$. This is the case because $(1+x)^ke^{-tx}$ is a
bounded function for $x\ge 0$ if $t>0$. Especially $s_t|_{\boundary
  M}=0$. (Here we use the fact that $\Delta$ is self adjoint on a
manifold of bounded geometry (
\ref{HsHodge})). Let $\psi:M\to [0,1]$ be a smooth function with
support contained in $X$ so that $\psi=1$ on an
eventually large set $Y\subset X$  and $\abs{\nabla\psi}<\epsilon$
(take f.i.\ an approximation of the Lipschitz function $x\mapsto
\max\{1,\epsilon \dist(x,M-X)\}$).
Observe that we find $M$-universal constants $C_{k,l}>0$ so that for $t>1$
\begin{equation}\label{Dst}
\begin{split}
 &\abs{\Delta^ls_t}_{H^k}\stackrel{\ref{HsHodge}}{=}\abs{(1+\Delta)^k\Delta^l s_t}_{L^2}\le C_{k,l} t^{-l}\abs{s}_{L^2},\\
 &\abs{\Delta^l s_t}_{C^k_b}\le C_{k,l} t^{-l}\abs{s}_{L^2}.
\end{split}\end{equation}
The first estimate follows from the spectral theorem, and the second one from the first and the Sobolev embedding theorem.\\
We arrive at the estimate
\begin{align*}
\int_{ Y}\abs{\nabla s_t}^2 &\le \int_{ M}\psi(\nabla s_t,\nabla s_t) = -(\nabla s_t,s_t\tensor\nabla\psi)_{L^2} + (\nabla s_t,\nabla(\psi s_t))_{L^2}\\
&= -(\nabla s_t,s_t\tensor\nabla\psi)_{L^2} + \underbrace{(\mathcal R^W s_t,\psi s_t)_{L^2}}_{\le 0} +\underbrace{\int_{\boundary M}(\mathcal S s_t,\psi s_t)}_{\le 0}\\
&\qquad + (d s_t,d(\psi s_t))_{L^2} + (\delta s_t,\delta(\psi s_t))_{L^2}\\
\intertext{(Apply \eqref{eq_Bochner}. Note $(\psi s_t)|_{\boundary M}=0$. Moreover, $\mathcal R^W\le 0$ and $\mathcal S\le 0$ on $\supp\psi$. For the next inequality use Cauchy-Schwartz and $\abs{\nabla\psi}\le\epsilon$:)}
 \le &\abs{\nabla s_t}_{L^2}\cdot\epsilon\cdot\abs{s_t}_{L^2}+ \abs{(ds_t,d\psi\wedge s_t)_{L^2}} + \abs{ds_t}^2_{L^2} + \abs{\delta s_t}^2_{L^2}\\
 & + \abs{(\delta s_t,*(d\psi\wedge *s_t))_{L^2}}\\
\stackrel{C_{1,0}}{\le}& \epsilon \abs{s}^2_{L^2} + \underbrace{\abs{ds_t}}_{\le\norm{d}_{\boundedops(H^1,L^2)}\abs{s_t}_{H^1}}\abs{s_t}_{L^2} \sup_{x\in M}{\underbrace{\abs{d\psi(x)}}_{=\abs{\nabla\psi}\le\epsilon}} \\
&+(\Delta s_t,s_t)_{L^2} + \abs{\delta s_t}_{L^2}\abs{s_t}_{L^2}\cdot\epsilon\\
\intertext{(First $\abs{\nabla s_t}\le \abs{s_t}_{H^1}\le C_{1,0}\abs{s}_{L^2}$ by \eqref{Dst}, then Cauchy-Schwartz and $(\Delta u,u)=(du,du)+(\delta u,\delta u)$ for $u\in\domain D_{\Delta}$ by \eqref{positivity})}
\le & \epsilon \abs{s}^2_{L^2} + \abs{\Delta s_t}_{L^2}\abs{s_t}_{L^2} +\epsilon\cdot(\norm{d}+\norm{\delta})C_{1,0}\abs{s}^2_{L^2}\quad\text{(by \eqref{Dst})}\\
\stackrel{C}{\le}  &\epsilon\abs{s}^2_{L^2} + t^{-1}\abs{s}_{L^2}^2
\qquad\text{(again by \eqref{Dst})} .
\end{align*}
The constant $C$ in the estimates is an $M$-universal constant involving $\norm{d}_{\boundedops(H^1,L^2)}$, $\norm{\delta}_{\boundedops(H^1,L^2)}$ and $C_{k,l}$ of \eqref{Dst} for $k,l\le 1$.\\
From now on, we proceed exactly as Roe does in \cite[1.11]{Roe(1988c)} to conclude
\begin{equation}\label{kl} \abs{s_t(y)}\le \epsilon_1(\epsilon)\abs{s}_{L^2} \quad\forall y\in Y_\epsilon 
\end{equation}
for $t$ sufficiently large,
where $\epsilon_1(\epsilon)\stackrel{\epsilon\to 0}{\to} 0$ and
$Y_\epsilon$ is an eventually large set depending only on
$\epsilon$. His analysis depends only on the following facts: 
\begin{itemize}
\item the estimates \eqref{Dst}
\item the existence of a function $V_0(r)$ with $\vol(B(x,r))\ge V_0(r)$ $\forall x\in M$, with $V_0$ monotonous and $V_0(r)\to\infty$ if $r\to\infty$.
\end{itemize}
Such a function can be constructed as follows:

Bounded geometry implies that we can construct $V_0(r)$ for $r$
sufficiently small \cite[Lemma
5]{Schick(1998b)}. To extend it to all of $\R_{\ge 0}$ we only have to show that 
\[ w(r):=\inf_{x\in M}\{\vol(B(x,r))\}\to\infty\quad\text{if $r\to\infty$} .\]
Actually, $w(r)$ grows at least linearly: Choose any $d>0$. Fix $r>0$ and $0<r_\epsilon<d/3$ so small that $B(x,r_\epsilon)$ is contained in some normal chart $\forall x\in M$. Then, we find $M$-universal $C>0$ so that $\vol(B(x,r_\epsilon))>C$ $\forall x\in M$. To given $x\in M$ choose $y\in M$ with $d(x,y)\ge 2r$. (this is possible because $M$ is complete but not compact). Choose a path from $x$ to $y$ with length $d(x,y)$. Mark successive points $x_0=x,x_1,\dots,x_{[r/d]}$ on the path with $d(x_i,x_{i+1})=d$.

Claim: $B(x_i,r_\epsilon)\cap B(x_j,r_\epsilon)=\emptyset$ if $i\ne j$. Else $d(x_i,x_j)<2r_\epsilon<d$, and the given path would not be length minimizing (in contradiction to its choice). Therefore 
\[\vol(B(x,r))\ge \sum_{i=0}^{[r/d]-1}\vol(B(x_i,r_\epsilon))\ge (r/d-1)\cdot C. \]

Having constructed all ingredients to use Roe's method, 
it remains to note that
$(1+\Delta)^ke^{-t\Delta}s\stackrel{t\to\infty}{\to}(1+\Delta)^kPs=Ps$
in $L^2$-norm $\forall k$.  Theorem \ref{HsHodge} shows that
$s_t\stackrel{t\to\infty}{\to}Ps$ in the Sobolev space $H^k$ $\forall
k$. The Sobolev embedding theorem  and  \eqref{kl}  imply 
\[ \abs{Ps(y)}\le \epsilon_1(\epsilon)\abs{s}_{L^2}\quad\forall y\in
Y_\epsilon . \quad\qed\]
\renewcommand{\qed}{}
\end{proof}

We will also need the following lemma in the proof of Theorem \ref{str_Boch}:
\begin{lemma}\label{cocompact}
  Let $M$ be a metric space with a cocompact action by isometries
  $\Gamma$. If $X\subset
  M$ is eventually large and $K\subset M$ is
  compact then there is $\gamma\in\Gamma$ with $\gamma(K)\subset X$.
\end{lemma}
\begin{proof}
  Choose $x\in K$. Then $\{B(x,n)\}_{n\in\N}$ is an open covering of
  $M$. Since the action of $\Gamma$ is by homeomorphism, the image of
  this covering under the projection $M\to M/\Gamma$ gives an open
  covering of the compact space $M/\Gamma$. Therefore we find $R>0$ so
  that every point in $M$ has a translate under $\Gamma$ which lies in
  $B(x,R)$. Let $D$ be the diameter of $K$. $D<\infty$ because $K$ is
  compact. Since $X$ is eventually large, there is $x_{D+R+1}$ with
  $B(x_{D+R+1},D+R+1)\subset X$. We just observed that
  $\gamma\in\Gamma$ exists with $\gamma(x)\subset B(x_{D+R+1},R)$
  ($\gamma$ is an isometry!). Then 
  \begin{equation*}
    \gamma(K)\subset
    \gamma(B(x,D))\subset B(x_{D+R+1},D+R+1)\subset X.\quad\qed
  \end{equation*}
\renewcommand{\qed}{}
\end{proof}

\begin{proof}[Proof of Theorem \ref{str_Boch}]
Let $D>0$ be a constant so that 
\[\begin{split}
& D^{-1}\abs{\omega}_{L^2(g_1)}\le \abs{\omega}_{L^2(g)}\le D\abs{\omega}_{L^2(g_1)}\\
\text{and } & D^{-1}\abs{\omega(x)}_{\Lambda^*(g_1)}\le\abs{\omega(x)}_{\Lambda^*(g)}\le D\abs{\omega(x)}_{\Lambda^*(g_1)} .
\end{split}\]
Note that $L^2( M, g)$ and $L^2( M,g_1)$ are equal as topological vector spaces since $ g$ and $g_1$ are bilipschitz.  The proof of Lemma \ref{pr_iso} shows that the restricted orthogonal projections 
\[\begin{split} P_1|: &L^2( M,g_1)\supset \mathcal H^p( M,\boundary M, g)\to \mathcal H^p( M,\boundary M, g_1)\\
P_2|: &L^2( M, g)\supset \mathcal H^p( M,\boundary M,g_1)\to \mathcal H^p( M,\boundary M, g)
\end{split}\]
are bounded and inverse to each other.\\
Suppose now $0\ne s\in\mathcal H^p( M,\boundary M, g)$ exists. Take $\abs{s}_{L^2(g)}=1$ and fix $x_0\in  M$ with $\abs{s(x_0)}=s_0>0$.
Let $p(x,y)$ be the (smooth) integral kernel of the projector $P_2$. Note that $Q(y)=p(x_0,y)$ is a smooth $L^2$-section of the bundle $\{\Hom(\Lambda^pT^*_y M,\Lambda^p T^*_{x_0} M)\}$ over $ M$. Choose $K$ compact so that 
\[\int_K\norm{Q(y)}^2\,dy\ge (1-\epsilon)\abs{Q}^2_{L^2}\quad\text{with}\quad 0<\sqrt{\epsilon}< \frac{s_0}{2\abs{Q}_{L^2}\norm{P_1}}.\]
Because $K$ is compact, $\vol(K)<\infty$. Choose 
\[0<\delta<\frac{s_0}{(2\abs{Q}_{L^2}D^2\sqrt{\vol(K)})}.\]
By \ref{roe} we find an eventually large set $Y\subset  M$ so that
\[ \abs{P_1f(x)}\le\delta\abs{f}_{L^2}\quad\forall x\in  Y\quad\forall f\in L^2 .\]
Since the isometry group $\Gamma$ acts cocompactly, we find by Lemma \ref{cocompact}
$\gamma\in\Gamma$ so that $\gamma^{-1}(K)\subset Y$.\\
Replace now $s$ by $\gamma^*s$, $x_0$ by $\gamma^{-1}(x_0)$, $Q(y)$ by $p(\gamma^{-1}x_0,y)$ and $K$ by $\gamma^{-1}K$. Because $\gamma$ acts isometric, for the new data holds 
\[s\in\mathcal H^p( M,\boundary M, g),\quad \abs{s(x_0)}=s_0,\qquad \int_K\norm{Q(y)}^2dy\ge (1-\epsilon)\abs{Q}^2_{L^2}.\]
 But now, by construction $K\subset Y$. Then (in terms of the metric $g$),
\[\begin{split} s_0^2 =& \abs{s(x_0)}^2= \abs{P_2(P_1)s(x_0)}^2 =\abs{\int_{ M} p(x_0,y)(P_1s)(y)\,dy}^2\\
\le &2\left(\abs{\int_K Q(y)P_1s(y)\,dy}^2+\abs{\int_{ M-K}Q(y)P_1s(y)\,dy}^2\right)\\
\le & 2\int_K\norm{Q(y)}^2\,dy \int_K\underbrace{\abs{P_1s(y)}^2}_{\le\delta D\abs{s}_{L^2(g_1)}}\,dy\\
& + 2\int_{ M-K}\norm{Q(y)^2\,dy} \underbrace{\int_{ M-K}{\abs{P_1s(y)}^2\,dy}}_{\le\abs{P_1s}^/_{L^2}}\\
\le  &2\abs{Q}^2_{L^2} \int_K D^2\delta^2\underbrace{\abs{s}^2_{L^2(g_1)}}_{\le D^2} +2\epsilon\abs{Q}^2_{L^2}\norm{P_1}^2\\
\le & 2\abs{Q}^2_{L^2}\delta^2 D^4\vol(K)
+2\epsilon\abs{Q}^2_{L^2}\norm{P_1}^2 <
\frac{s_0^2}{2}+\frac{s_0^2}{2}= s_0^2 .
\end{split}\]
This is a contradiction, therefore no non-trivial relative $L^2$-harmonic forms exist.\\
For the absolute case, simply replace $\mathcal S$ by $\mathcal S*$ from equation \eqref{eq_Bochner}.
\end{proof}

\begin{corollary}
Let $M$ be a compact manifold (possibly with boundary). Let $\tilde M$
be an orientable infinite normal covering of $M$. In the following, we
are considering only metrics of bounded geometry.
\begin{itemize}
\item If the first relative $L^2$-Betti number $b^1_{(2)}(\tilde
  M,\boundary \tilde M)\ne 0$ then there is no metric $g_1$ on
  $\tilde M$ which lies in the natural bilipschitz class, so that
  $\Ric\ge 0$  and $\tr l\ge 0$ on any eventually large sets. In
particular this is the case if $\chi(M)<0$ and either $\dim M=2$ or $\dim
M=4$ and $\boundary M=\emptyset$.
\item If the relative Euler characteristic $\chi(M,\boundary M)\ne 0$
  then in the natural bilipschitz class on the universal covering
  $\tilde M$ no metric exists so that the curvature operator and the
  second fundamental form operator are both $\ge 0$ on eventually large sets.
\end{itemize}
\end{corollary}
\begin{proof}
These are direct consequences of the algebraic results in Proposition
\ref{pos_def} as soon as we note that because of the $L^2$-index
theorem for coverings of $\boundary$-manifolds $\chi(M)=\sum_p (-1)^p
b^p_{(2)}(\tilde M)$ and $\chi(M,\boundary M)=\sum_p (-1)^p
b^p_{(2)}(\tilde M,\boundary \tilde M)$ (compare \cite[Proposition
6.4]{Schick(1998c)}). Moreover, if $M$ is even-dimensional then
because of Poincare duality $\chi(M)=\chi(M,\boundary M)$. The
statement concerning closed $4$-manifold also follows from Poincare
duality which implies $b^1_{(2)}(\tilde M)=b^{3}_{(2)}(\tilde M)$.
\end{proof}

\begin{remark} 
John Roe \cite{Roe(1988c)} proves a  theorem which implies the same for
 coverings of manifolds without boundary and amenable covering group. 
Roe does not use eventually large, but more
exhaustive sets --- in particular complements of compact sets. His
question whether the theorem holds for arbitrary coverings, is
answered affirmatively by our result.
\end{remark}


%% file: L2_DodziukSchick.tex
\section{Proof of the $L^2$-de Rham theorem}\label{sec_Rham}

The proof of Theorem \ref{Hodge_de_Rham_theorem} is done in three steps.

\subsubsection*{Step 1: $\boundary M=\emptyset$}
This is a classical result of Dodziuk \cite{Dodziuk(1981)}, compare
also \cite{Dodziuk(1977)}.

\subsubsection*{Step 2: product metrics}\label{sec_prod}
Let $(M,g)$ be an oriented manifold of bounded geometry with boundary $\boundary M=M_1\amalg M_2$. Suppose the Riemannian metric is a product near the boundary. Construct the following ``quadruple":
\[ W= (M\cup_{M_1} -M)\cup_{M_2\amalg -M_2}-(M\cup_{M_1}-M) .\]
We have an obvious $V_4=\Z/2\times \Z/2$-action on $W$ generated by flips $\tau_1$ along $M_1\amalg -M_1$ and $\tau_2$ along $M_2\amalg-M_2$.

The Riemannian metric on $M$ can be extended to a metric $\bar g$ on
$W$ so that $V_4$ acts isometrically. This induces a unitary $V_4$
action on $\mathcal{H}(W)$. Use the $g$-bounded triangulation $K$ of
$M$ to get a $V_4$-invariant triangulation $\bar K$ of $W$, which is $\bar g$-bounded. Then  we get a unitary action $V_4$ on $H^{(*)}_{(2)}( W)$. Since these operations come from diffeomorphism of the underlying manifold the de Rham map $A$ is $V_4$-equivariant by its geometric definition. Dodziuk's theorem says that $A$ is an isomorphism
\[ A:\mathcal{H}^*( W)\to H^*_{(2)}(\bar K) .\]

Restriction to $V_4$-invariant subspaces  yields
that the de Rham map restricted to the $(-+)$-eigenspaces is an isomorphism
\[ A: \mathcal{H}^{-+}( W) \to H^{-+}_{(2)}(\bar K) .\]
Here $X^{-+}$ is characterized by $\tau_1 x=-x$ and $\tau_2 x=x$ for $x\in X$.

It remains to identify these eigenspaces.
\begin{lemma}
Denote the inclusion of one copy of $ M$ in $ W$ by $i: M\hookrightarrow W$. Then we get an isomorphism
\[ i^*: \mathcal{H}^{-+}( W) \to \mathcal{H}( M, M_1) .\]
\end{lemma}
\begin{proof}
We simply write down the inverse map: take $\omega\in \mathcal{H}( M,
M_1)$ $\implies \omega|_{ M_1}=0=(*\omega)|_{ M_2}$. Define
$\bar\omega\in L^2( W)$ by $i^*\bar\omega=\omega$,
$i^*\tau_1^*\bar\omega=-\omega$, $i^*\tau_2^*\bar\omega=\omega$,
$i^*\tau_1^*\tau_2^*\bar\omega=-\omega$. Obviously, $V_4$ acts in the
correct way on $\bar\omega$. We have to check that $\bar\omega$ is not
only in $L^2$ but actually in $\mathcal{H}$. It suffices to show that
$\Delta\bar\omega=0$ in the weak sense. Then by elliptic regularity
(\ref{elliptic_regularity} and \ref{unif_Laplace}) $\bar\omega\in
H^\infty$ and $\Delta\bar\omega=0$ smoothly. That weakly $\Delta\bar\omega=0$ is shown using integration by parts.

\end{proof}

\begin{lemma}The inclusion $i: M\hookrightarrow  W$ induces an inclusion $i:K\hookrightarrow\bar K$ and this an isomorphism
\[ i^*: H^{-+}_{(2)}(\bar K) \to H_{(2)}(K, K_1) . \]
\end{lemma}
\begin{proof}
The existence of an inverse map on the level of $L^2$-cochains is, given the boundary conditions, clear. It induces the inverse map on cohomology.
\end{proof}

\subsubsection*{Step 3: general case}

We have the following strategy: If the metric $g$ is not a product
near the boundary, the doubling trick does not work since the metric
would no longer be smooth. However, we now construct a metric $\bar g$
which is a product near the boundary and which allows  an isomorphism
of the spaces of harmonic forms commuting with the de Rham map.  Therefore, the positive result for product metrics implies the general result.

\begin{proposition}\label{p_met}
Let $(M,g)$ be a manifold of bounded geometry with boundary $\boundary M$. Suppose $K: \boundary M\times [0,r_C)\to M$ is a metric collar. Let $g'=K_*(g|_{\boundary M}\times g_{\R})$ be the corresponding product metric on the collar.\\
Set $\bar g:=\varphi g'+ (1-\varphi)g$, where $\varphi:[0,r_C)\to [0,1]$ is $1$ in a neighborhood of $0$ and $0$ outside another neighborhood of $0$.\\
Then $(M,\bar g)$ is a manifold of bounded geometry, every bundle
$E\bundleover M$ which is a bundle of bounded geometry with respect to
$g$ has bounded geometry also with respect to $\bar g$, and the
Sobolev spaces and their topology do not depend on the metric: $H^s(E,g)=H^s(E,\bar g)$.
\end{proposition}
\begin{proof}
First, observe $g|_{\boundary M}=\bar g|_{\boundary M}$. Moreover, the
inward unit normal vector fields $\nu$ and $\bar\nu$ are equal. The
geodesics in $(M,g)$ with starting velocity $\nu$ and length $<r_C$
are by the very construction of $g'$ also geodesics of
$(\im(K),g')$. Let $\alpha(t)$ be such a geodesic with
$v_t=\alpha'(t)$. Then $\alpha$ is a geodesic of $(M,\bar g)$ too,
since for arbitrary vector fields $Z$ (compare Gallot/Hulin/Lafontaine
\cite[2.51]{Gallot-Hulin-Lafontaine(1987)}):
\[\begin{split}
& 2\bar g(\nabla^{\bar g}_{v_t}v_t,Z) = 2 v_t.\bar g(v_t,Z)-Z.\bar g(v_t,v_t)-2\bar g([v_t,Z],v_t)\\
&\quad=
\varphi 2g'(\underbrace{\nabla^{g'}_{v_t}v_t}_{=0},Z) + (1-\varphi)2g(\underbrace{\nabla^g_{v_t}v_t}_{=0},Z)
+ 2v_t(\varphi)(g'(v_t,Z)\\ 
&\quad-g(v_t,Z)) +Z(\varphi)(g'(v_t,v_t)-g(v_t,v_t)) = 0 .
\end{split}\]
The latter is true since the Gauss lemma 
 says $g(v_t,Z)=0\iff g'(v_t,Z)=0$ and moreover
$g(v_t,v_t)=1=g'(v_t,v_t)$. This is proved in \cite[2.93]{Gallot-Hulin-Lafontaine(1987)}
for radial geodesics, and the prove with not more than the obvious
changes applies also to our product situation.
Apparently $\nabla^{\bar g}_{v_t}v_t=0$ and $\alpha$ is a
geodesic. \\
It follows that normal boundary coordinates for $(M,g)$ and $(M,\bar
g)$ coincide. If we take $\supp\varphi$ and the injectivity constant
$r_i$ to be sufficiently small, the same holds for Gaussian
coordinates centered at a point in $M-N_{1/3}$ (because $g$ and $\bar
g$ are equal there).

To show that $(M,\bar g)$ is a manifold of bounded geometry, it
suffices to compute the metric tensors $\bar g_{ij}$ and $\bar g^{ij}$
in (the unchanged!) normal boundary coordinates and find $M$-universal
bounds for them and their partial derivatives. For $\bar
g_{ij}(x',t)=\varphi(t) g_{ij}(x',0)+(1-\varphi(t))g_{ij}(x',t)$ this
is immediate. Now the product rule for derivatives implies that it
suffices to find bounds for $\abs{\bar g^{ij}}$ only. Equivalently,
produce a bound for the norm of the matrix $(\bar g^{ij})_{ij}$, i.e.\
a lower bound for the smallest eigenvalue of the positive self adjoint
matrix $(\bar g_{ij})_{ij}$. By assumption, we find corresponding
$M$-universal bounds for $(g_{ij}(x',t))_{ij}$. The smallest
eigenvalue of a positive self adjoint matrix $G$ is equal to
$\inf_{\abs{x}=1} (Gx,x)$. For convex combinations, if $(G_1x,x)\ge c$
and $(G_2x,x)\ge c$ then also $((\varphi G_1+(1-\varphi)G_2)x,x)\ge
c$.\\ 
This implies that we get an $M$-universal bound for $\abs{(\bar
g^{ij})}$; and $(M,\bar g)$ has bounded geometry.\\ 
For the definition of the Sobolev spaces $H^s(E)$, we have seen now
that we can choose identical data for $g$ and for $\bar g$.
\end{proof}


To compare the spaces of harmonic forms (and the de Rham map) for two different metrics we use the following functional analytical lemma:
\begin{lemma}\label{Banach_projection}
Suppose $X$ is a Banach space, $U$, $V$ and $A$ are closed subspaces with
\[ U\cap A=\{0\}=V\cap A; \qquad U+A=X=V+A .\]
Then the projection $p_U:X=U+A \to U$ onto $U$ along $A$, when
restricted to $V$ is an isomorphism $V\to U$ of topological vector spaces with inverse $p_V$.
\end{lemma}

Now suppose that $(M,g_1)$ is a complete oriented  Riemannian $\boundary$-manifold and $g_2$ is another metric on $M$ bilipschitz to $g_1$. Operators and spaces  depending on the metric will be decorated with subscripts $1$ or $2$ accordingly. Note that the topological vector space $L^2(\Lambda^p(T^* M))$ does not depend on the metric.

\begin{lemma}\label{pr_iso}  The $g_2$-orthogonal projection
\[ pr_2: L^2(\Lambda^p(T^* M))\to \mathcal{H}^p( M, M_1; g_2) \]
restricted to $\mathcal{H}( M, M_1; g_1)$ is a bounded isomorphism with inverse the corresponding restriction $pr_1|$.
\end{lemma}
\begin{proof}
We want to apply Lemma \ref{Banach_projection}. Therefore, we prove \[\mathcal{H}^p( M, M_1;g_1)\oplus\overline{d\Omega_{d}^{p-1}} = \mathcal{H}^p( M, M_1;g_2)\oplus\overline{d\Omega_{d}^{p-1}}.\]
This is true because it is exactly the kernel of $d$, which does not depend on 
the metric. The same is true for 
$\overline{ d C_0^\infty(\Lambda^{p-1}(T^* M);\;\cdot|_{ M_1}=0)}=\overline{ d\Omega_d^{p-1}}$, the image of $d$.
\end{proof}

\begin{lemma} Suppose we are in the situation of Lemma \ref{pr_iso} and, additionally, $(M,g_1)$ and $(M,g_2)$ are manifolds of bounded geometry such that all the Sobolev spaces do not depend on the metric. Let $K$ be a $g_1$-bounded triangulation of $M$. Then the projections $pr_i|$ commute with the de Rham map $A$.
\end{lemma}
\begin{proof}
Write $h\in \mathcal{H}^p( M, M_1; g_1)$ as $h=h'+x$ with $h'\in\mathcal{H}^p( M, M_1; g_2)$ and $x\in \overline{d\Omega_{\infty}^{p-1}}$ (use the Hodge decomposition \ref{HsHodge}). $h,h'\in H^\infty$ $\implies x\in H^\infty$ and $x|_{ M_1}=0$ since $h|_{ M_1}=0=h'|_{ M_1}$. Moreover, $dx=0$. Therefore, the following statement concludes the proof.
\end{proof}
\begin{lemma}\label{pr_dR}
Suppose $x\in H^\infty\cap \overline{d\Omega_{\infty}^{p-1}}$ and $x|_{ M_1}=0$ and $dx=0$. Then $[Ax]=0\in H^p_{(2)}( M, M_1)$, where $A$ is the de Rham map.
\end{lemma}
\begin{proof}
The de Rham map 
\[ A: H^{2k}(\Lambda^p(T^* M)) \to C^p_{(2)}(K) \]
is defined and bounded for $k$ sufficiently large. 
(Here $C^p_{(2)}$ denotes the space of square summable cochains of the  triangulation $K$. The proof of Dodziuk \cite{Dodziuk(1977),Dodziuk(1981)} remains valid).
It induces a bounded map
\[ A:\{\omega\in H^{2k}(\Lambda^p(T^* M)); \omega|_{
M_1}=0,\,d\omega=0\} \to H^p_{(2)}( M, M_1) .\]
Since $x\in H^\infty$, $Ax\in H^p_{(2)}( M, M_1)$ is defined.

By Theorem \ref{HsHodge}, $x=(1+\Delta)^k y$ for a unique $y\in \overline{d\Omega^{p-1}_\infty}$ (from elliptic regularity we know that actually $y\in H^\infty$). Remember that here the closure is taken in the topology of $H^{2k}$.  We find 
\[y_n\in H^\infty(\Lambda^{p+1}(T^* M), bd(\infty))\quad\text{with}\quad y=\lim_n dy_n.\]
 Especially $dy =\lim_n ddy_n = 0$.
Since $(1+\Delta)^k=\sum_{i=0}^k\binom{k}{i}\Delta^i$ and 
\begin{gather*} \Delta^i =(d\delta+\delta d)^i \stackrel{dd=0=\delta\delta}{=}(d\delta)^i+(\delta d)^i\quad\text{for }i>0 \\
\implies x=(1+\Delta)^k y= y+\sum_{i=1}^k
\binom{k}{i}\left((d\delta)^iy+\underbrace{(\delta d)^i y}_{=0}\right)
= y +dz  \end{gather*}
with $z=\delta\sum_{i=1}^k \binom{k}{i} (d\delta)^{i-1} y$.
Note that $z\in H^\infty$ since $y\in H^\infty$. $H^{2k}$-boundedness
of the boundary differential
operators implies $z|_{M_1}=0$ as $y_n$ fulfill 
$bd(\infty)$. This implies $A(dz)=d(Az)=0\in H^p_{(2)}( M, M_1)$
\[\begin{split} \implies Ax &=Ay+A(dz)=Ay.\qquad\text{ Now}\\
 Ay &=A(\lim_n dy_n)\stackrel{A \text{ is $H^{2k}$-continuous}}{=}\lim_n A(dy_n)
 \end{split}\]
and $A(dy_n)=dA(y_n)=0\in H^p_{(2)}( M, M_1)$.
\end{proof}

The metric $\bar g$ constructed in Lemma \ref{p_met} fulfills all the conditions of the lemma and has product structure near the boundary. For this by Step 2 the de Rham map is an isomorphism, and therefore the Hodge-De Rham theorem is proved in general.

\begin{remark}
Dodziuk \cite{Dodziuk(1981)} uses $L^2$-de Rham cohomology to prove his theorem. This can also be introduced for $\boundary$-manifolds of
bounded geometry. 
There are  several $L^2$-de Rham complexes whose cohomology is isomorphic 
to the simplicial cohomology  defined above, as shown in 
\cite[Chapter 6]{Schick(1996)}.
\end{remark}

